\documentclass[11pt,twoside]{article}

\usepackage{latexsym}
\usepackage{amsmath}
\usepackage{amsthm}
\usepackage{amssymb}
\usepackage{vmargin}
\usepackage{amscd}
\usepackage{stmaryrd}
\usepackage{euscript}
\usepackage{mathrsfs}
\usepackage{amscd}
\usepackage[all]{xy}
\usepackage{xr}
\DeclareMathAlphabet{\mathpzc}{OT1}{pzc}{m}{it}

\setmargins{32mm}{20mm}{14.6cm}{22cm}{1cm}{1cm}{1cm}{1cm}

\newtheorem{theorem}{Theorem}[section]
\newtheorem{proposition}[theorem]{Proposition}
\newtheorem{corollary}[theorem]{Corollary}

\newtheorem{lemma}[theorem]{Lemma}
\newtheorem*{theorem*}{Theorem}
\newtheorem*{proposition*}{Proposition}
\newtheorem*{corollary*}{Corollary}
\newtheorem*{lemma*}{Lemma}
\newtheorem*{conjecture*}{Conjecture}

\theoremstyle{definition}
\newtheorem{definition}[theorem]{Definition}

\newtheorem*{definition*}{Definition}

\theoremstyle{remark}

\newtheorem{remark}[theorem]{Remark}
\newtheorem{remarks}[theorem]{Remarks}

\newtheorem*{example*}{Example}
\newtheorem*{examples*}{Examples}
\newtheorem*{remark*}{Remark}
\newtheorem*{remarks*}{Remarks}
\newtheorem*{exercise*}{Exercise}

\newcommand\da{\!\downarrow\!}

\newcommand\id{\mathrm{id}}

\newcommand\ten{\otimes}

\newcommand\CC{\mathrm{C}}

\newcommand\RR{\mathrm{R}}

\newcommand\Ru{\mathrm{R_u}}

\renewcommand\H{\mathrm{H}}
\newcommand\z{\mathrm{Z}}

\newcommand\Z{\mathbb{Z}}

\newcommand\R{\mathbb{R}}
\newcommand\Cx{\mathbb{C}}

\newcommand\vv{\mathbb{V}}
\newcommand\ww{\mathbb{W}}

\newcommand\bB{\mathbb{B}}

\newcommand\bG{\mathbb{G}}

\newcommand\bL{\mathbb{L}}

\newcommand\bO{\mathbb{O}}

\newcommand\bS{\mathbb{S}}

\newcommand\cE{\mathcal{E}}

\newcommand\cH{\mathcal{H}}

\newcommand\cM{\mathcal{M}}
\newcommand\cN{\mathcal{N}}

\renewcommand\O{\mathscr{O}}

\newcommand\sA{\mathscr{A}}

\newcommand\sE{\mathscr{E}}
\newcommand\sF{\mathscr{F}}

\newcommand\sO{\mathscr{O}}

\newcommand\m{\mathfrak{m}}

\newcommand\g{\mathfrak{g}}

\newcommand\fh{\mathfrak{h}}

\newcommand\Ho{\mathrm{Ho}}

\newcommand\Alg{\mathrm{Alg}}

\newcommand\Hom{\mathrm{Hom}}

\newcommand\End{\mathrm{End}}

\newcommand\Aut{\mathrm{Aut}}
\newcommand\Out{\mathrm{Out}}
\newcommand\ROut{\mathrm{ROut}}

\newcommand\VHS{\mathrm{VHS}}

\newcommand\Gal{\mathrm{Gal}}

\newcommand\im{\mathrm{Im\,}}

\newcommand\Ob{\mathrm{Ob}\,}

\newcommand\Top{\mathrm{Top}}

\newcommand\agp{\mathrm{AGp}}
\newcommand\agpd{\mathrm{AGpd}}

\newcommand\mal{\mathrm{Mal}}

\newcommand\Spec{\mathrm{Spec}\,}

\newcommand\Aff{\mathrm{Aff}}

\newcommand\sch{\mathrm{sch}}

\newcommand\Sing{\mathrm{Sing}}

\newcommand\ad{\mathrm{ad}}
\newcommand\norm{\mathrm{norm}}

\newcommand\Lim{\varprojlim}

\newcommand\into{\hookrightarrow}

\newcommand\xra{\xrightarrow}

\newcommand\alg{\mathrm{alg}}

\newcommand\dmd{\diamond}
\newcommand\bt{\bullet}
\newcommand\by{\times}

\newcommand\mc{\mathrm{MC}}

\newcommand\Gg{\mathrm{Gg}}

\newcommand\Rep{\mathrm{Rep}}

\newcommand\GL{\mathrm{GL}}

\newcommand\hol{\mathrm{hol}}

\newcommand\diag{\mathrm{diag}\,}

\newcommand\pd{\partial}
\newcommand\dc{d^{\mathrm{c}}}
\newcommand\half{\frac{1}{2}}

\newcommand\gpd{\mathrm{Gpd}}

\renewcommand\alg{\mathrm{alg}}
\newcommand\red{\mathrm{red}}

\newcommand\co{\colon\thinspace}

\renewcommand\Alg{\mathrm{Alg}}

\begin{document}
\title{Non-abelian real Hodge theory for proper varieties}
\author{J.P.Pridham\thanks{The author is supported by Trinity College, Cambridge}}
\maketitle

\begin{abstract}
 We show that if $X$ is any proper complex variety, there is a weight decomposition on the real schematic homotopy type, in the form of an algebraic $\bG_m$-action. This extends to a real Hodge structure, in the form of a discrete $\Cx^*$-action, such that $\Cx^* \by X \to X^{\sch}$ is real analytic. If the fundamental group is algebraically good, and the higher homotopy groups have finite rank, this gives bigraded decompositions on the complexified homotopy groups. For smooth proper varieties, the Hodge structure can be recovered from the cohomology ring with coefficients in the universal semisimple local system. 
\end{abstract}

\tableofcontents

\section*{Introduction}
\addcontentsline{toc}{section}{Introduction}

In \cite{Simpson}, Simpson established the existence of a pure Hodge structure on the real reductive pro-algebraic fundamental group $\varpi_1(X,x)^{\red}$ of a compact K\"ahler manifold $X$. This took the form of a discrete action of the circle group $U(1)$ on $\varpi_1(X,x)^{\red}$, in such a way that the resulting map
$$
U(1) \by \pi_1(X,x) \to \varpi_1(X,x)^{\red}(\R)
$$
was continuous. In fact, the results of  \cite{Sim2} imply that this map is  real holomorphic (i.e. corresponds to a $\Gal(\Cx/\R)$-equivariant holomorphic map on the $\Cx$-valued points of these schemes).

The purpose of this paper is to establish the existence of a real Hodge structure on the pro-algebraic homotopy types as defined in \cite{htpy} (or equivalently, on the schematic homotopy types of \cite{chaff}) of complex varieties. The definition and main properties of pro-algebraic homotopy types are summarised in Section \ref{review}.  

In \cite{Hodge2} Definition 2.1.4, Deligne defined a real Hodge structure on a real vector space to be an algebraic action of the real algebraic group $\Cx^*$. A pure Hodge structure of weight $0$ is then a $\Cx^*$-action for which $\R^* < \Cx^*$ acts trivially, giving an action of $\Cx^*/\R^* \cong U(1)$.
Accordingly, we look to extend Simpson's $U(1)$-action on $\varpi_1(X,x)^{\red}(\R)$ to a $\Cx^*$-action on the whole real pro-algebraic homotopy type $G(X)^{\alg}$ (the algebraisation of the path groupoid of $X$). 
Our notion of real real Hodge structure will be a  real holomorphic map (or, rather, a homotopy class of maps)
$$
\Cx^* \by G(X) \to G(X)^{\alg}(\R),
$$
such that the induced map
$$
\Cx^* \to \Hom(G(X), G(X)^{\alg}(\R)) \cong \End(G(X)^{\alg})
$$
is multiplicative. We also require that the reductive pro-algebraic fundamental groupoid $\varpi_f(X)^{\red}$ is of weight zero under this action, and that the induced map
$$
\R^* \by G(X)^{\alg} \to G(X)^{\alg}
$$   
is algebraic (i.e. a morphism of schemes). In Section \ref{RHS}, this definition  is made precise, and  several properties of non-abelian Hodge structures are given. In particular, if the group $\pi_1(X,x)$ is algebraically good, and the higher homotopy groups have finite rank, this gives algebraic $\Cx^*$-actions on the real vector spaces $\pi_n(X,x)\ten_{\Z}\R$ (and hence bigraded decompositions on the $\pi_n(X,x)\ten_{\Z}\Cx$).

In Section \ref{cK}, we establish the existence of a real  Hodge structure on the homotopy type of a compact K\"ahler manifold. This is done by combining the $U(1)$-action of \cite{Simpson} with Hodge theory for semisimple local systems. Formality of the pro-algebraic homotopy type (from the $d\dc$-lemma)  means that it is determined by the cohomology ring of the universal semisimple local system. Characterisation of cohomology groups by harmonic forms then enables to construct  a $\Cx^*$-action on cohomology, by studying behaviour of the Laplacian under Simpson's $U(1)$-action. However, we have not been able to compare this Hodge structure fully with the  Hodge structures of \cite{KTP} and \cite{Morgan}.

Section \ref{propcx} is concerned with extending the non-abelian real Hodge structures from smooth proper varieties to singular proper varieties, in the spirit of Hodge III (\cite{Hodge3}). The key idea is that we may replace any proper variety by a smooth proper simplicial variety. The machinery of \cite{htpy} is well-suited to dealing with simplicial spaces and cohomological descent, so functoriality of the non-abelian Hodge structure on smooth proper varieties allows it to extend naturally to all proper varieties.

\section{Review of pro-algebraic homotopy types}\label{review}
Here we give a summary of the results from \cite{htpy} which will be needed in this paper. Fix a field $k$ of characteristic zero.

\subsection{Pro-algebraic groupoids}
We first recall some definitions from \cite{htpy} \S\S \ref{htpy-gpdsn}--\ref{htpy-levisn}.
\begin{definition}
Define a pro-algebraic groupoid $G$ over $k$ to consist of the following data:
\begin{enumerate}
\item A discrete set $\Ob(G)$.
\item For all $x,y \in \Ob(G)$, an affine scheme $G(x,y)$ (possibly empty) over $k$.
\item A groupoid structure on $G$, consisting of an associative multiplication  morphism $m:G(x,y)\by G(y,z) \to G(x,z)$, identities $\Spec k \to G(x,x)$ and inverses $G(x,y) \to G(y,x)$
 \end{enumerate}
Note that a pro-algebraic group is just a pro-algebraic groupoid on one object.
We say that a pro-algebraic groupoid is reductive (resp. pro-unipotent) if the pro-algebraic groups $G(x,x)$ are so for all $x \in \Ob(G)$. An algebraic groupoid is a pro-algebraic groupoid for which the $G(x,y)$ are all of finite type.
\end{definition}
If $G$ is a pro-algebraic groupoid, let $O(G(x,y))$ denote the global sections of the structure sheaf of $G(x,y)$.

\begin{definition}
Given morphisms $f,g:G \to H$ of pro-algebraic groupoids, define a natural isomorphism $\eta$ between $f$ and $g$ to consist of morphisms
$$
\eta_x: \Spec k \to H(f(x),g(x))
$$
for all $ x\in \Ob(G)$, such that the following diagram commutes, for all $x,y \in \Ob(G)$:
$$
\begin{CD}
G(x,y) @>f(x,y)>> H(f(x),f(y))\\
@Vg(x,y)VV  @VV{\cdot\eta_y}V \\
 H(g(x),g(y)) @>{\eta_x\cdot}>>  H(f(x),g(y)).
\end{CD}
$$
 
A morphism $f:G \to H$ of pro-algebraic groupoids is said to be an equivalence if there exists a morphism $g:H \to G$ such that $fg$ and $gf$ are both naturally isomorphic to identity morphisms. This is the same as saying that for all $y \in \Ob(H)$, there exists $x \in \Ob(G)$ such that $H(f(x),y)(k)$ is non-empty (essential surjectivity), and that for all $x_1,x_2 \in \Ob(G)$, $G(x,y) \to G(f(x_1),f(x_2) )$ is an isomorphism.
\end{definition}

\begin{definition}\label{gpdrep}
Given a pro-algebraic groupoid $G$, define a finite-dimensional linear $G$-representation to be a functor $\rho:G \to \mathrm{FDVect}_k$ respecting the algebraic structure. Explicitly, this  consists of a set $\{V_x\}_{x \in \Ob(G)}$ of finite-dimensional $k$-vector spaces, together with morphisms $\rho_{xy}:G(x,y) \to \Hom(V_y,V_x)$ of affine schemes, respecting the multiplication and identities. 

A morphism $f:(V,\rho)\to (W,\varrho)$ of $G$-representations consists of $f_x \in \Hom(V_x,W_x)$ such that 
$$
f_x\circ\varrho_{xy}=\rho_{xy}\circ f_y:G(x,y) \to \Hom(V_x,W_y).
$$
\end{definition}

\begin{definition}
Given a pro-algebraic groupoid $G$, define the reductive quotient $G^{\red}$ of $G$ by setting $\Ob(G^{\red})=\Ob(G)$, and
$$
G^{\red}(x,y)=G(x,y)/\Ru(G(y,y))= \Ru(G(x,x))\backslash G(x,y),
$$
where $\Ru(G(x,x))$ is the pro-unipotent radical of the pro-algebraic group $G(x,x)$.
The equality arises since if $f\in G(x,y),\, g \in \Ru(G(y,y))$, then $f gf^{-1}\in \Ru(G(x,x))$, so both equivalence relations are the same. Multiplication and inversion descend similarly. Observe that $G^{\red}$ is then a reductive pro-algebraic groupoid. Representations of $G^{\red}$ correspond to semisimple representations of $G$.
\end{definition}

\begin{definition}
Let $\agpd$ denote the category of pro-algebraic groupoids over $k$, and observe that this category is contains all (inverse) limits. There is functor from $\agpd$ to $\gpd$, the category of abstract groupoids, given by $G \mapsto G(k)$. This functor preserves all limits, so has a left adjoint, the algebraisation functor, denoted $\Gamma \mapsto \Gamma^{\alg}$. This can be given explicitly by $\Ob(\Gamma)^{\alg}=\Ob(\Gamma)$,  and 
$$
\Gamma^{\alg}(x,y)=\Gamma(x,x)^{\alg}\by^{\Gamma(x,x)}\Gamma(x,y),
$$
where $\Gamma(x,x)^{\alg}$ is the pro-algebraic completion of the group $\Gamma(x,x)$. 

The finite-dimensional linear representations of $\Gamma$ (as in Definition \ref{gpdrep}) correspond to those of $\Gamma^{\alg}$, and these can be used to recover $\Gamma^{\alg}$, by Tannakian duality.
\end{definition}

\begin{definition}
Given a pro-algebraic groupoid $G$, and $U=\{U_x\}_{x \in \Ob(G)}$ a collection of pro-algebraic groups parametrised by $\Ob(G)$, we say that $G$ acts on $U$ if there are morphisms $ U_x\by G(x,y) \xra{*} U_y$ of affine schemes, satisfying the following conditions:
\begin{enumerate}
\item $(uv)*g= (u*g)(v*g)$, $1*g=1$ and $(u^{-1})*g= (u*g)^{-1}$,  for $g \in G(x,y)$ and $u,v \in U_x$.

\item $u*(gh)=(u*g)*h$ and  $u*1=u$, for $g \in G(x,y), h \in G(y,z)$ and $u \in U_x$. 
\end{enumerate}

If $G$ acts on $U$, we write $G \ltimes U $ for the groupoid given by  
\begin{enumerate}
\item $\Ob(G \ltimes U):=\Ob(G)$.

\item $( G \ltimes U)(x,y):= G(x,y)\by U_y$.

\item $(g,u)(h,v):= (gh, (u*h)v)$ for $g \in G(x,y), h \in G(y,z)$ and $u \in U_y, v \in U_z$. 
\end{enumerate}
\end{definition}

\begin{definition}
Given a pro-algebraic groupoid $G$, define $\Ru(G)$ to be the collection $\Ru(G)_x=\Ru(G(x,x))$ of pro-unipotent pro-algebraic groups, for $x \in \Ob(G)$. $G$ then acts on $\Ru(G)$ by conjugation, i.e.
$$
u*g:= g^{-1}u g,
$$
for $u \in \Ru(G)_x$, $g \in G(x,y)$.
\end{definition}

\begin{proposition}
For any pro-algebraic groupoid $G$, there is a Levi decomposition $G=G^{\red} \ltimes \Ru(G)$, unique up to conjugation by $\Ru(G)$.
\end{proposition}
\begin{proof}
\cite{htpy} Proposition \ref{htpy-leviprop}.
\end{proof}

\subsection{The pro-algebraic homotopy type of a topological space}
We now recall the results from \cite{htpy} \S \ref{htpy-sagpdsn}.

\begin{definition}
Let $\bS$ be the category of simplicial sets, and $s\gpd$ the category of simplicial groupoids on a constant set of objects (as in \cite{sht}). Let $\Top$ denote the category of compactly generated Hausdorff topological spaces. 

A map $f:X \to Y$ in  $\Top$ is said to be a weak equivalence if it gives an isomorphism $\pi_0X \to \pi_0Y$ on path components, and for all $x \in X$, the maps $\pi_n(f):\pi_n(X,x) \to \pi_n(Y,fx)$ are all isomorphisms. A map $f:X \to Y$ in $\bS$ is said to be a weak equivalence if the map $|f|:|X| \to |Y|$ is so. A map $f:G \to H$ in $s\gpd$ is a weak equivalence if the map on components $\pi_0G_0 \to \pi_0H_0$ is an isomorphism, and for all objects $x \in \Ob G$, the maps $\pi_n(G(x,x))\to \pi_n(H(x,x))$ are all isomorphisms.  

For each of these categories, we define the corresponding homotopy categories $\Ho(\bS),\Ho(s\gpd),\Ho(\Top)$ by localising at weak equivalences.
\end{definition}

Note that there is a functor from $\Top$ to $\bS$ which sends $X$ to the simplicial set
$$
\Sing(X)_n= \Hom_{\Top}(|\Delta^n|, X).
$$
this gives an  equivalence of the corresponding homotopy categories, whose quasi-inverse is geometric realisation. From now on, we will thus restrict our attention to simplicial sets.

As in \cite{sht} Ch.V.7, there is a classifying space  functor $\bar{W}:s\gpd \to \bS$, with left adjoint $G:\bS \to s\gpd$, Dwyer and Kan's  path groupoid functor (\cite{pathgpd}), and these give  equivalences $\Ho(\bS)\sim \Ho(s\gpd)$.  The geometric realisation of $|G(X)|$ is weakly equivalent to the path space of $|X|$. These functors have the additional properties that $\Ob G(X)=X_0$, $(\bar{W}G)_0=\Ob(G)$, $\pi_0G(X)\cong \pi_0|X|$, $\pi_0(|\bar{W}G|)\cong \pi_0G_0$,   $\pi_n(G(X)(x,x))\cong\pi_{n+1}(|X|,x)$ and $\pi_{n+1}(|\bar{W}G|,x)=\pi_n(G(x,x))$. This allows us to study simplicial groupoids instead of topological spaces.

\begin{definition}
Given a simplicial object $G_{\bullet}$ in the category of pro-algebraic groupoids,   with $\Ob(G_{\bt})$ constant, define the fundamental groupoid $\pi_f(G_{\bullet})$ of $G_{\bt}$ to have objects $\Ob(G)$, and for $x,y \in \Ob(G)$, set 
$$
\pi_f(G)(x,y):= G_0(x,y) / \sim,
$$ 
where $\sim$ is the equivalence relation generated by $\pd_0h\sim \pd_1h$ for $h \in G(x,y)$. This is also pro-algebraic.
\end{definition}

\begin{definition}
Define a pro-algebraic simplicial  groupoid  to consist of a simplicial complex $G_{\bullet}$ of pro-algebraic groupoids, such that $\Ob(G_{\bt})$ is constant  and for all $x \in \Ob(G)$, $G(x,x)_{\bt} \in s\agp$, i.e. the maps $G_n(x,x) \to \pi_0(G)(x,x)$ are pro-unipotent extensions of pro-algebraic groups.  We denote the category of pro-algebraic simplicial groupoids by $s\agpd$. 

Define a morphism $f:G_{\bullet} \to H_{\bullet}$ in $s\agpd$ to be a weak equivalence if the map  $\pi_f(f):\pi_f(G_{\bullet})\to \pi_f(H_{\bullet})$ is an equivalence of pro-algebraic groupoids, and the maps  $\pi_n(f,x):\pi_n(G_{\bullet}(x,x)) \to \pi_n(H_{\bullet}(fx,fx))$ are isomorphisms for all $n$ and for all $x \in \Ob(G)$. We define $\Ho(s\agpd)$ to be the localisation of $s\agpd$ at weak equivalences.
\end{definition}

There is a forgetful functor $(k):s\agpd \to s\gpd$, given by sending $G_{\bt}$ to $G_{\bt}(k)$. This functor  has a left adjoint $G_{\bt} \mapsto (G_{\bt})^{\alg}$. We can describe $(G_{\bt})^{\alg}$ explicitly. First let $(\pi_f(G))^{\alg}$ be the pro-algebraic completion of the abstract groupoid $\pi_f(G)$, then let $(G^{\alg})_n$ be the relative Malcev completion  (defined in \cite{malcev} for pro-algebraic groups) of the morphism
$$
G_n \to (\pi_f(G))^{\alg}.
$$
In other words, $G_n \to (G^{\alg})_n \xra{f} (\pi_f(G))^{\alg}$ is the universal diagram with $f$ a pro-unipotent extension.

\begin{proposition}\label{algqd}
The functors $(k)$ and ${}^{\alg}$ give rise to a pair of adjoint functors 
$$
\xymatrix@1{\Ho(s\gpd) \ar@<1ex>[r]^{\bL^{\alg}} & \Ho(s\agpd) \ar@<1ex>[l]^{(k)}_{\bot} },
$$
with $\bL^{\alg}G(X)=G(X)^{\alg}$, for any $X \in \bS$.
\end{proposition} 
\begin{proof}
\cite{htpy} Proposition \ref{htpy-algqd}.
\end{proof}

\begin{definition}
Given a  simplicial set (or equivalently a  topological space), define the pro-algebraic homotopy type of $X$ over $k$ to be the object
$$
G(X)^{\alg}
$$ 
in $\Ho(s\agpd)$. Define the pro-algebraic fundamental groupoid by $\varpi_f(X):=\pi_f(G(X)^{\alg})$. Note that  $\pi_f(G^{\alg})$ is the pro-algebraic completion of the fundamental groupoid $\pi_f(G)$.

We then define the higher homotopy groups $\varpi_n(X)$  (as $\varpi_fX$-representations) by 
$$
\varpi_n(X):=\pi_{n-1}(G(X)^{\alg}),
$$
where $\pi_n(G)$ is the representation $x \mapsto \pi_n(G(x,x))$, for $x \in \Ob(G)$.
\end{definition}

\subsection{Relative Malcev homotopy types}\label{malcev}

\begin{definition}\label{malcevdef}
Assume we have an abstract groupoid $G$, a reductive pro-algebraic groupoid $R$, and a representation $\rho:G \to R(k)$ which is an isomorphism on objects and Zariski-dense on morphisms (i.e. $\rho\co G(x,y) \to R(k)(\rho x, \rho y)$ is Zariski-dense for all $x,y \in \Ob G$). Define the Malcev completion $(G,\rho)^{\mal}$ (or $G^{\rho, \mal}$)   of $G$ relative to $\rho$  to be the universal diagram
$$
G \to (G,\rho)^{\mal} \xra{p} R,
$$
with $p$  a pro-unipotent extension, and the composition equal to $\rho$. Explicitly, $\Ob(G,\rho)^{\mal}=\Ob G$ and 
$$
(G,\rho)^{\mal}(x,y)=(G(x,x), \rho)^{\mal}\by^{G(x,x)}G(x,y).
$$
If $G$ and  $R$ are groups, observe that this agrees with the usual definition.  

If $\varrho:G \to R(k)$  is any any  Zariski-dense representation (i.e. essentially surjective on objects and Zariski-dense on morphisms) to a reductive pro-algebraic groupoid (in most examples, we take $R$ to be a group), we can define another reductive groupoid $\tilde{R}$ by setting $\Ob \tilde{R}=\Ob G$, and $ \tilde{R}(x,y)=R(\varrho x, \varrho y)$. This gives a representation $\rho:\pi_fX \xra{\rho} \tilde{R}$ satisfying the above hypotheses, and we define the Malcev completion of $G$ relative to $\varrho$ to be the Malcev completion of $G$ relative to $\rho$. Note that $\tilde{R} \to R$ is an equivalence of pro-algebraic groupoids.
\end{definition}

\begin{definition}
Given a Zariski-dense morphism $\rho:\pi_fX \to R(k)$, let the Malcev completion $G(X,\rho)^{\mal}$ of $X$ relative to $\rho$ be the pro-algebraic simplicial group $(G(X), \rho)^{\mal}$.  Observe that the Malcev completion of $X$ relative to $(\pi_fX)^{\red}$ is just $G(X)^{\alg}$. Let $\varpi_f(X,\rho)^{\mal}=\pi_fG(X,\rho)^{\mal}$ and $\varpi_n(X,\rho)^{\mal}=\pi_{n-1}G(X,\rho)^{\mal}$. Note that $\pi_f((X,\rho)^{\mal})$ is the relative Malcev completion of $\pi_f\rho:\pi_fX \to R(k)$.  
\end{definition} 

\begin{definition}
Define a groupoid $\Gamma$ to be good with respect to a Zariski-dense representation $\rho: \Gamma \to R(k)$ to a reductive pro-algebraic groupoid if the map
$$
\H^n(\Gamma^{\rho, \mal}, V) \to \H^n(\Gamma, V)
$$  
is an isomorphism for all $n$ and all finite-dimensional $\Gamma^{\rho, \mal}$-representations $\Gamma$.
\end{definition}

\begin{lemma}
Assume that for all $x \in \in \Ob \Gamma$,  $\Gamma(x,x)$ is finitely presented, with $\H^n(\Gamma, - )$ commuting with filtered direct limits of  $\Gamma^{\rho, \mal}$-representations, and $\H^n(\Gamma, V)$ finite-dimensional for all finite-dimensional $\Gamma^{\rho, \mal}$-representations $V$. 

Then $\Gamma$ is good with respect to $\rho$ if and only if for any finite-dimensional $\Gamma^{\rho, \mal}$-representation $V$, and $\alpha \in \H^n(\Gamma, V)$, there exists an injection $f:V \to W_{\alpha}$ of finite-dimensional $\Gamma^{\rho, \mal}$-representations, with $f(\alpha)=0 \in \H^n(\Gamma, W_{\alpha})$.
\end{lemma}
\begin{proof}
As for \cite{schematic} Lemma 4.15.
\end{proof}

\begin{theorem}\label{classicalpimal}
If $X$ is a  topological space with fundamental groupoid $\Gamma$, equipped with a Zariski-dense representation $\rho: \Gamma \to R(k)$ to a reductive pro-algebraic groupoid  for which: 
\begin{enumerate}
\item $\Gamma$ is  good with respect to $\rho$,
\item $\pi_n(X,-)$ is of finite rank for all $n>1$,
\item and the $\Gamma$-representation  $\pi_n(X,-)\ten_{\Z} k$ is an extension of $R$-representations (i.e. a $\Gamma^{\rho, \mal}$-representation),
\end{enumerate}
then the canonical map
$$
  \pi_n(X,-)\ten_{\Z} k \to \varpi_{n}(X^{\rho, \mal},-) 
$$
is an isomorphism for all $n>1$.
\end{theorem}
\begin{proof}
\cite{htpy} Theorem \ref{htpy-classicalpimal}.
\end{proof}

\subsection{Equivalent formulations}

Fix a reductive pro-algebraic groupoid $R$.

\subsubsection{Simplicial Lie algebras}

\begin{definition}
Given a $k$-algebra $A$,  let $\hat{\cN}_A$ be opposite to the category of ind-conilpotent Lie  coalgebras, and let  $\hat{\cN}_A(R)$ be the category of $R$-representations in $\hat{\cN}_A$.Write $s\hat{\cN}_A(R)$ for the category of simplicial objects in $\hat{\cN}_A(R)$. A weak equivalence in $s\hat{\cN}_A(R)$ is a map which gives isomorphisms on cohomology groups of the duals  (which are just $A$-modules). We denote by $\Ho(s\hat{\cN}_A(R))$ the localisation of  $s\hat{\cN}_A(R)$ at weak equivalences. For $k=A$, we will usually drop the subscript, so $\hat{\cN}(R):=\hat{\cN}_k(R)$, and so on.
\end{definition}

\begin{definition}
Define $\cE(R)$ to be the full subcategory of $\agpd\da R$ consisting of those morphisms $\rho:G\to  R$ of  proalgebraic groupoids which are pro-unipotent extensions. Similarly, define $s\cE(R)$ to consist of the pro-unipotent extensions in $s\agpd\da R$, and $\Ho(s\cE(R))$ to be the localisation of $s\cE(R)$ at weak equivalences.
\end{definition}

\begin{definition}
Given a pro-algebraic groupoid $R$, define the category $s\cM_A(R)$ to have the same objects as $s\hat{\cN}_A(R)$, with morphisms given by
$$
\Hom_{s\cM_A(R)}(\g,\fh)=\Hom_{\Ho(s\hat{\cN}_A(R))}(\g,\fh)/\exp(\fh_0^R),
$$
where $\fh_0^R$ is the sub-Lie algebra of $\fh_0(A)$ fixed by $R$, acting by conjugation on the set of homomorphisms.
\end{definition}

\begin{proposition}\label{meequiv}
For any reductive pro-algebraic groupoid $R$, the categories $\Ho(s\cE(R))$ and $s\cM(R)$ are equivalent, sending $\g \in s\hat{\cN}(R)$ to $R\ltimes \exp(\g)$.
\end{proposition}
\begin{proof}
\cite{htpy} Proposition \ref{meequiv}. 
\end{proof}

\begin{definition}
We can now define the Malcev homotopy type of $X$ relative to $\rho$ to be the image of   $G(X,\rho)^{\mal}$ in $\Ho(s\cE(\tilde{R}))$, or equivalently $\Ru G(X,\rho)^{\mal}$ in  $s\cM(\tilde{R})$. Since $\tilde{R} \to R$ is an equivalence of groupoids, there is an equivalence $s\hat{\cN}(R) \to s\hat{\cN}(\tilde{R})$, so may regard the Malcev homotopy type as an object of $s\cM(R)$ (or of $\Ho(s\cE(R))$). It pro-represents the functor
$$
\g \mapsto \Hom_{\Ho(\bS\da BR(k))}(X,\bar{W}(R(k)\ltimes \exp(\g))), 
$$
for $\g \in s\cM(R)$.
\end{definition}

\subsubsection{Chain Lie algebras}

\begin{definition}Let  $dg\hat{\cN}_A$  be opposite to the category of non-negatively graded ind-conilpotent cochain Lie  coalgebras over $A$. Define $dg\hat{\cN}_A(R)$ to be the category of $R$-representations in $dg\hat{\cN}_A$.
A weak equivalence in $dg\hat{\cN}_A(R)$ is a map which induces isomorphisms on cohomology groups of the duals. We denote by $\Ho(dg\hat{\cN}_A(R))$ the localisation of  $dg\hat{\cN}_A(R)$ at weak equivalences. For $k=A$, we will usually drop the subscript, so $dg\hat{\cN}(R):=dg\hat{\cN}_k(R)$, and so on.
\end{definition}

\begin{definition}
We say that a morphism $f:\g \to \fh$ in $dg\hat{\cN}(R)$ is free if there exists a (pro-finite-dimensional)  sub-$R$-representation $V \subset \fh$ such that  $\fh$ is the free pro-nilpotent graded Lie algebra over $\g$ on generators $V$.
\end{definition}

\begin{proposition}[Minimal models]\label{dgminimal}
For every object $\g$ of   $dg\hat{\cN}(R)$, there exists a free chain Lie algebra $\m$  with $d=0$ on  the abelianisation $\m/[\m,\m]$,  unique up to non-unique isomorphism, together with a weak equivalence  $\m \to \g$.
\end{proposition}
\begin{proof}
 \cite{htpy} Proposition \ref{htpy-dgminimal}. 
\end{proof}

\begin{definition}
Let $dg\cM_A(R)$ be the category with the same objects as $dg\hat{\cN}_A(R)$, and morphisms given by
$$
\Hom_{dg\cM_A(R)}(\g,\fh)= \Hom_{\Ho(dg\hat{\cN}_A(R))}(\g,\fh)/\exp(\fh^R_0).
$$
\end{definition}

\begin{proposition}\label{nequiv}
There is a normalisation functor $N:s\hat{\cN}_A(R) \to dg\hat{\cN}_A(R)$ such that
$$
\H_i(N\g) \cong \pi_i(\g),
$$
giving equivalences $\Ho(s\hat{\cN}_A(R))\simeq \Ho(dg\hat{\cN}_A(R))$, and $s\cM_A(R) \simeq dg\cM_A(R)$.
\end{proposition}
\begin{proof}
\cite{htpy} Propositions \ref{htpy-nequiv} and \ref{htpy-anequiv}.
\end{proof}

\subsubsection{Cosimplicial algebras}

\begin{definition}
 Let $c\Alg(R)$ be the category of of $R$-representations in cosimplicial $k$-algebras. A weak equivalence in $c\Alg(R)$ is a map which induces isomorphisms on cohomology groups. We denote by $\Ho(c\Alg(R))$ the localisation of  $c\Alg(R)$ at weak equivalences.
\end{definition}

\begin{definition}
Let $s\Aff(R)$ denote the category  of simplicial affine schemes over $k$, i.e. the category opposite to $c\Alg(R)$. Similarly, let $\Ho(s\Aff(R))$ be the category opposite to $\Ho(c\Alg(R))$.
\end{definition}

\begin{definition}
Given $V,W \in \Rep(R)$, define $V\ten^RW:=\Hom_{\Rep(R)}(k,V\ten W)$.
\end{definition}

\begin{definition} Given $A \in c\Alg(R)$ and $\g \in s\hat{\cN}(R)$, define the 
 Maurer-Cartan space $\mc(A,G)$ to consist of sets $\{\omega_n\}_{n\ge 0}$, with $\omega_n \in \exp(A^{n+1}\hat{\ten}^R\g_n)$, such that 
\begin{eqnarray*}
\pd_i\omega_n &=& \left\{\begin{matrix} \pd^{i+1}\omega_{n-1}  & i>0 \\ (\pd^1\omega_{n-1})\cdot(\pd^0\omega_{n-1})^{-1} & i=0,\end{matrix} \right.\\
\sigma_i\omega_n &=& \sigma^{i+1}\omega_{n+1},\\
\sigma^0\omega_n&=& 1,
\end{eqnarray*}
where  $\exp(A^{n+1}\hat{\ten}^R\g_n)$ is the group  whose underlying set is the Lie algebra $A^{n+1}\hat{\ten}\g_{n-1}$, with multiplication given by the Campbell-Baker-Hausdorff formula.
\end{definition}

\begin{definition}
Given $A \in c\Alg(R)$ and $\g \in s\hat{\cN}(R)$, define the gauge group $\Gg(A,\g)\le \prod_n \exp(A^n\hat{\ten}^R\g_n)$
to consist of those $g$ satisfying
\begin{eqnarray*}
\pd_ig_n &=& \pd^{i}g_{n-1}  \quad \forall i>0, \\
\sigma_ig_n &=& \sigma^{i}g_{n+1} \quad \forall i.
\end{eqnarray*}
 This has an action on $\mc(A,\g)$ given by  
$$
(g*\omega)_n= (\pd_0g_{n+1}) \cdot \omega_n \cdot (\pd^0g_n^{-1}),
$$
and we define the torsor space by
$$
\pi(A,\g):=\mc(A,\g)/\Gg(A,\g).
$$
\end{definition}

\begin{definition}
Let $c\Alg(R)_0$ be the full subcategory of $c\Alg(R)$ whose objects satisfy $\H^0(A) \cong k$. Let $\Ho(c\Alg(R)_0)$ be the full subcategory of $\Ho(c\Alg(R)_0)$ with objects in $c\Alg(R)_0$. Let $s\Aff(R)_0$ be the category opposite to $c\Alg(R)_0$, and   $\Ho(s\Aff(R)_0)$ opposite to $\Ho(c\Alg(R)_0)$.
\end{definition}

\begin{proposition}\label{wequiv}
There is a pair of equivalences
$$
\xymatrix@1{ \Ho(s\Aff(R))_0  \ar@<1ex>[r]^-{\bar{G}}  & s\cM(R) \ar@<1ex>[l]^-{\bar{W}}},
$$
given by
$$
\Hom_{\Ho(s\Aff(R))}(\Spec A,\bar{W}\g)=\Hom_{s\cM(R)}(\bar{G}(A), \g)=\pi(A,\g).
$$
\end{proposition}
\begin{proof}
\cite{htpy} Proposition \ref{htpy-wequiv}.
\end{proof}

\begin{definition}
Given a topological space $X$, and a sheaf $\sF$ on $X$, define 
$$
\CC^n(X,\sF):= \prod_{f:\Delta^n \to X} \Gamma(\Delta^n, f^{-1}\sF).
$$
Together, these form a cosimplicial complex $\CC^{\bt}(X,\sF)$.
\end{definition}

\begin{definition}
Recall that  $O(R)$ has the natural structure of an $R\by R$-representation. Since  every $R$-representation has an associated semisimple local system on $|BR(k)|$, we will also write $O(R)$ for the $R$-representation in semisimple local systems on $|BR(k)|$  corresponding to the $R\by R$-representation $O(R)$.
We then define the $R$-representation $\bO(R)$ in semisimple local systems on $X$ by $\bO(R):=\rho^{-1}O(R)$.
\end{definition}

\begin{proposition}
Under the equivalences of Propositions \ref{meequiv} and \ref{wequiv}, the relative Malcev homotopy type $G(X)^{\rho, \mal}$ of a topological space $X$ corresponds to 
$$
\CC^{\bt}(X,\bO(R))\in c\Alg(R).
$$
\end{proposition}\label{eqhtpy}
\begin{proof}
\cite{htpy} Theorem \ref{htpy-eqhtpy}
\end{proof}

\begin{corollary}\label{eqtoen}
Pro-algebraic homotopy types are equivalent to the schematic homotopy types of \cite{chaff}, in the sense that the full subcategory of the homotopy category $\Ho(s\mathrm{Pr})$ on objects $X^{\sch}$ is equivalent to the full subcategory of $\Ho(s\agpd)$ on objects $G(X)^{\alg}$. Under this equivalence, pro-algebraic homotopy groups are isomorphic to schematic homotopy groups.
\end{corollary}
\begin{proof}
\cite{htpy} Corollary \ref{htpy-eqtoen}.
\end{proof}

\subsubsection{Cochain algebras}

\begin{definition}
Define $DG\Alg(R)$ to be the category of $R$-representations in  non-negatively graded cochain $k$-algebras. 
A weak equivalence in $DG\Alg(R)$ is a map which induces isomorphisms on cohomology groups. We denote by $\Ho(DG\Alg(R))$ the localisation of  $DG\Alg(R)$ at weak equivalences.
Define $dg\Aff(R)$ to be the category opposite to $DG\Alg(R)$,  and $\Ho(dg\Aff(R))$  opposite to $\Ho(DG\Alg(R))$. 

Let $DG\Alg(R)_0$ be the full subcategory of $DG\Alg(R)$ whose objects $A$ satisfy $\H^0(A)=k$. Let $\Ho(DG\Alg(R))_0$ be the full subcategory of $\Ho(DG\Alg(R))$ on the objects of  $DG\Alg(R)_0$. Let $dg\Aff(R)_0$ and $\Ho(dg\Aff(R))_0$ be the opposite categories to $DG\Alg(R)_0$ and $\Ho(DG\Alg(R))_0$, respectively.
\end{definition}

\begin{proposition}\label{affequiv}
There is a denormalisation functor $D:DG\Alg(R) \to c\Alg(R)$ such that
$$
\H^i(DA) \cong \H^i(A),
$$
giving an equivalence  $\Ho(c\Alg(R))\simeq \Ho(DG\Alg(R))$.
\end{proposition}
\begin{proof}
\cite{htpy} \ref{htpy-affequiv}.
\end{proof}

\begin{definition}
Given a cochain algebra $A \in DG\Alg(R)$, and a chain Lie algebra $\g \in dg\hat{\cN}(R)$, define 
$$
\mc(A,\g):=\{\omega \in \bigoplus_n A^{n+1}\hat{\ten}^R \g_n \,|\,d\omega+\half[\omega,\omega]=0\}.
$$
\end{definition}

\begin{definition}\label{dgdef}\label{dgdefgauge}
Given $A \in DG\Alg(R)$ and $\g \in dg\hat{\cN}(R)$, we define the gauge group by
$$
\Gg(A,\g):= \exp(\prod_n A^n\hat{\ten}^R\g_n).
$$ 
Define a gauge action of $\Gg(A,\g)$ on $\mc(A,\g)$ by 
$$
g(\omega):= g\cdot \omega \cdot g^{-1} -(dg)\cdot g^{-1}.
$$

We define the
torsor space by
$$
\pi(A,\g)=\mc(A,\g)/\Gg(A,\g).
$$
\end{definition}

\begin{theorem}\label{bigequiv}
We have the following commutative diagram of equivalences of categories:
$$
\xymatrix{
\Ho(dg\Aff(R))_0 \ar@<1ex>[r]^{\Spec D} \ar@<1ex>[d]^{\bar{G}}& Ho(s\Aff(R))_0  \ar@<1ex>[d]^{\bar{G}} \\
dg\cM(R) \ar@<1ex>[u]^{\bar{W}}  & s\cM(R), \ar@<1ex>[u]^{\bar{W}} \ar@<1ex>[l]^{N},
}
$$
with the pair
$$
\xymatrix@1{ \Ho(dg\Aff(R))_0  \ar@<1ex>[r]^-{\bar{G}}  & dg\cM(R) \ar@<1ex>[l]^-{\bar{W}}},
$$
given by
$$
\Hom_{\Ho(dg\Aff(R))}(\Spec A,\bar{W}\g)=\Hom_{dg\cM(R)}(\bar{G}(A), \g)=\pi(A,\g).
$$
\end{theorem}
\begin{proof}
\cite{htpy} Theorem \ref{htpy-qs}, Corollary \ref{htpy-bigequiv} and Theorem \ref{htpy-defqs}. 
\end{proof}

\begin{definition}
Given a manifold $X$, denote the sheaf  of real $n$-forms on $X$ by $\sA^n$. Given a real sheaf $\sF$ on $X$, write
$$
A^n(X,\sF):=\Gamma(X,\sF\ten_{\R} \sA^n).
$$ 
\end{definition}

\begin{proposition}\label{propforms}
The real Malcev  homotopy type of a manifold $X$ relative to $\rho:\pi_fX \to R(\R)$ is given in $DG\Alg(R)$ by 
$
A^{\bt}(X, \bO(R)).
$
\end{proposition}
\begin{proof}
\cite{htpy} Proposition \ref{htpy-propforms}.
\end{proof}

\subsection{Weights}\label{wgts}

\begin{definition}
Given $\g \in s\hat{\cN}_R$ and a $k$-algebra $A$, define 
$$
\Out_R(\g)(A):= \Aut_{s\cM_A(R)}(\g\hat{\ten} A).
$$

Given $G \in s\cE(R)$, define $\ROut(G):=\Out_R(\Ru(G))$, noting that $\ROut(G)(k) \cong \Aut_{\Ho(s\cE(R))}(G)$.
 For $G \in s\agpd$, set $\ROut(G):=\Out_{G^{\red}}(\Ru(G))$. 
\end{definition}

\begin{definition}
By \cite{htpy} Theorem \ref{htpy-auto}, if $\H^i(X,V)$ is finite-dimensional for all finite-dimensional irreducible $R$-representations $V$ and $\rho:G(X) \to R$, then
$$
\ROut(G(X)^{\rho,\mal})
$$
is represented by  a pro-algebraic group over $k$, and we define a weight decomposition on $G(X)^{\rho,\mal}$ to be a morphism
$$
\bG_m \to \ROut(G(X)^{\rho,\mal})
$$
of pro-algebraic groups.
\end{definition}

\begin{definition}
Given a reductive pro-algebraic groupoid $R$, we say that $A \in DG\Alg(R)$ is formal if it is weakly equivalent to its cohomology algebra $\H^*(A)$. We say that the Malcev homotopy type $(X,\rho)^{\mal}$ of a topological space $X$  relative to a Zariski-dense homomorphism $\rho:\pi_fX \to R$ is formal if  the corresponding DG algebra  (given by Theorem \ref{bigequiv}) is formal.
\end{definition}

\begin{proposition}\label{formalpin}
If $G(X,\rho)^{\mal}$ is formal, then it has a canonical weight decomposition, for which $\H^i(X,\vv)$ is pure of weight $i$.
\end{proposition}
\begin{proof}
\cite{htpy} Corollary \ref{htpy-formalpin}.
\end{proof}

\section{Real Hodge structures}\label{RHS}
From now on, we will take $k=\R$. 
We now wish  to study the Hodge structure on cohomology. 

\subsection{Review of classical real Hodge structures}

In this section, we recall the standard definitions of real Hodge theory, and fix some notation.

\begin{definition}
Recall from \cite{Hodge2} Definition 2.1.4 that a real  Hodge structure on a real vector space $V$ is an action of the real algebraic group $S:=\Cx^*$, obtained from $\bG_m$ by Weil restriction of scalars from $\Cx$ to $\R$. Explicitly, for any real algebra $A$,
$$
S(A)=\{(x,y)\in A^2\,|\, x^2+y^2 \text{ is invertible}\}.
$$
Therefore 
\begin{eqnarray*}
S(\Cx)&\cong& \Cx^*\by \Cx^*\\ 
(a,b) &\mapsto& (a+ib, a-ib).
\end{eqnarray*}
\end{definition}

For a vector space $V$, with complexification $V_{\Cx}$ this corresponds to a Hodge decomposition 
$$
V_{\Cx}=\bigoplus_{p+q}V_{\Cx}^{pq}
$$  
by setting
$$
\lambda'\mu''(\sum v_{pq})=\sum \lambda^p\mu^q v_{pq},
$$
for $(\lambda,\mu) \in  (\Cx^*)^2  \cong S(\Cx)$.

Since $\overline{V_{\Cx}^{pq}}=\overline{V_{\Cx}}^{qp}$ is the condition for this Hodge structure to be real, note that 
$$
\overline{\lambda'\mu''v}=\bar{\mu}'\bar{\lambda}''\bar{v},
$$
 so that this  does indeed descend to a real  $S$-action. 

\begin{definition}\label{dmd}
For $\lambda \in  \Cx^* \cong S(\R) $, we will denote the action of $\lambda$ on $V$ by
$$
\lambda \dmd v := \lambda'\bar{\lambda}''v.
$$
\end{definition}

Given a $\Cx$-linear map $F:V \to V$ of type $(a,b)$, observe that 
$$
\lambda \dmd Fv =\lambda^a \bar{\lambda}^b F(\lambda \dmd v).
$$

\subsection{Non-abelian real Hodge structures}

\begin{definition}\label{outchar}
For $G \in s\agpd$, define  the group $\Out(G)$ of outer automorphisms of $G$ to consist of pairs $\theta=(\theta^{\red},\Ru(\theta))$, with $\theta^{\red}:G^{\red} \to G^{\red}$ an automorphism, and 
$$
\Ru(\theta) \in \Hom_{\Ho(s\agpd\da G^{\red})}(\theta^{\red}_*G, G)
$$
an isomorphism, where $\theta^{\red}_*G$ is the composition $G \to G^{\red} \xra{\theta^{\red}} G^{\red}$.

We define the group structure on $\Out(G)$ by 
$$
\theta \circ \phi =(\theta^{\red} \circ \phi^{\red}, \Ru(\theta)\circ \theta^{\red}_*\Ru(\phi)). 
$$
\end{definition}

\begin{definition}
Given a real affine scheme $X$, define the ring $O(X)^{\hol}$ of real holomorphic functions on $X$ by
$$
O(X)^{\hol}:= \H^0(X(\Cx), \sO^{\hol}_{X(\Cx)})/\Gal(\Cx/\R),
$$
where $\sO^{\hol}_{X(\Cx)}$ is the sheaf of complex holomorphic functions on $X(\Cx)$.

Given real affine schemes $X, Y$, define a real holomorphic map $f:X \to Y$ to be a $\Gal(\Cx/\R)$-equivariant holomorphic map $f:X(\Cx) \to Y(\Cx)$. Observe that this is equivalent to an element of $Y(O(X)^{\hol})$.
\end{definition}

\begin{definition}
Given $X \in \bS$, $H \in s\agpd_{\R}$, and a $k$-algebra $A$, we define the set $O\Hom(X,\bar{W}H)(A)$  of $A$-valued outer homomorphisms to be the coequaliser
$$
\xymatrix@1{\Hom_{s\gpd}(X,\bar{W}PH(A)) \ar@<0.5ex>[r] \ar@<-0.5ex>[r]& \Hom_{\bS}(X,\bar{W}H(A)) \ar[r] & O\Hom(X,\bar{W}H)(A)}.
$$
Here,  $PH:=H^I\by_{(H^{\red})^I} H^{\red}$, for 
$H^I$  a path object for $H$ in $s\agpd$; equivalently $PH$ is a path object for $H$ in $s\agpd \da H^{\red}$.

Note that if $\rho:\pi_fX \to R$ is as in Definition \ref{malcevdef}, then 
$$
\Out(G(X)^{\rho,\mal})\subset O\Hom(X,\bar{W}G(X)^{\rho,\mal})(\R).
$$

For a real holomorphic space $Z$, we define the set $O\Hom(X,\bar{W}H)(Z)$  of $Z$-valued outer homomorphisms by
$$
O\Hom(X,\bar{W}H)(Z):= O\Hom(X,\bar{W}H)(O(Z)^{\hol}).
$$
\end{definition}

\begin{lemma}
$$
O\Hom(X,\bar{W}H)(A)= \mc(X,H(A))/\Gg(X,\Ru(H)(A)).
$$
\end{lemma}

\begin{definition}\label{realhodgedef}
Given $X \in \bS$ and a Zariski-dense homomorphism, isomorphic on objects, $\rho: \pi_fX \to H(\R)$   to a reductive pro-algebraic groupoid over $\R$,     such that $\H^i(X, \rho^*V)$ is finite-dimensional for all $i$ and all finite-dimensional irreducible $H$-representations $V$,  define a real Hodge structure on $\rho$ to consist of 
an $S$-valued outer homomorphism
$$
h \in O\Hom(X, \bar{W}G(X)^{\rho, \mal})(S_{\hol}),
$$ 
such that 
\begin{enumerate}

\item for all $s \in S(\R)$, the morphism $h(s)^{\red}:\varpi_f(X)^{\red} \to H$ descends to an automorphism of $H$, and    the induced map
$$
h^{\delta}:S(\R) \to     O\Hom(G(X)^{\rho, val}, G(X)^{\rho, \mal})
$$
is a monoid homomorphism, or equivalently
$$
h^{\delta}:S(\R)^{\delta} \to \Out(G(X)^{\rho, \mal}) 
$$
is a group homomorphism;

\item the action of $S(\R)$ on $(G(X)^{\rho, \mal})^{\red}=H$ is pure of weight $0$, i.e. the composition
$$
\bG_m(\R) \into S(\R) \to \Out(G(X)^{\rho, \mal}) \xra{p} \Aut(H)
$$ 
is trivial;

\item the map
$$
h:\bG_m \to \ROut(G(X)^{\rho, \mal})=\ker p
$$
is a homomorphism of pro-algebraic groups. In other words, the restriction of $h$ to $\bG_m$ is a weight decomposition on $G(X)^{\rho, \mal}$ in the sense of \S \ref{wgts}.
\end{enumerate}

Define a real Hodge structure on $X$ to be a real Hodge structure on the canonical map $\rho: \pi_fG \to \varpi_f(X)^{\red}(\R)$, in which case $G(X)^{\rho, \mal}=G(X)^{\alg}$. 
\end{definition}

\begin{remark}
 Note that, unlike the definition of pure Hodge structures in \cite{Simpson}, we do not have any hypothesis on Hodge type. This is because this would not be satisfied by singular varieties.
\end{remark} 

On passing to homotopy groups, we obtain the following.
\begin{proposition}
For any point $x$ in a space  $X$ with a real Hodge structure, there are real holomorphic maps
$$
S \by \pi_n(X,x) \to \varpi_n(X,x),
$$
unique up to conjugation by $\Ru(\varpi_1(X,x))$,  such that the induced map 
$$
S(\R) \to \End(\varpi_n(X,x)(\R))/\Ru(\varpi_1(X,x))
$$ 
is a monoid homomorphism, and such that the composition
$$
\bG_m \into S \to \Out(\varpi_n(X,x)(\R))
$$
is a homomorphism of pro-algebraic groups. In particular, this gives weight decomposition on the pro-algebraic homotopy groups, unique up to conjugation by $\Ru(\varpi_1(X,x))$. 
\end{proposition}

\subsection{Variations of Hodge structure}\label{varhodgesn}

The following definition is adapted  from \cite{weight1}:
\begin{definition}
Given  a discrete group $\Gamma$ acting on a pro-algebraic groupoid $G$, for which the action on $\Ob G$ is trivial,
define ${}^{\Gamma}\!G$ to be the maximal quotient of $G$ on which $\Gamma$ acts algebraically.  This is the inverse limit $\Lim_{\alpha} G_{\alpha}$ over those surjective maps
$$
G \to G_{\alpha},
$$
with $G_{\alpha}$ algebraic, for which the $\Gamma$-action descends to $G_{\alpha}$. Equivalently, $O( {}^{\Gamma}\!G)$ is the sum of those finite-dimensional $\Gamma$-representations of $O(G)$ which are closed under the coproduct.
\end{definition}

\begin{definition}
Given  a real Hodge structure on $\rho: \pi_fX \to \H(\R)$, define 
$$
{}^{\VHS}\!H:= {}^{S^{\delta}}\!H.
$$
\end{definition}

\begin{remarks}
This notion is analogous to the definition given in \cite{weight1} of the maximal quotient on which Frobenius acts algebraically. In the same way that representations of that group corresponded to semisimple subsystems of local systems underlying Weil sheaves, representations of ${}^{\VHS}\!\varpi_f(X)^{\red}$ will correspond to 
local systems underlying variations of Hodge structure. Also note that since $\bG_m(\R) < S(\R)$ acts trivially on $\varpi_f(X)^{\red}$, the $S(\R)$-action factors through the circle group $U(1)\cong S/\bG_m$. 
\end{remarks}

\begin{proposition}\label{vhsalg}
The action of $S/\bG_m$ on ${}^{\VHS}\!H$ is algebraic, in the sense that
$$
S/\bG_m \by {}^{\VHS}\!H\to {}^{\VHS}\!H
$$
is a morphism of schemes.

It is also an inner action, coming from a morphism
$$
S/\bG_m  \to (\prod_{x \in X_0} {}^{\VHS}\!H(x,x))/Z({}^{\VHS}\!H)
$$
of pro-algebraic groupoids, where $Z$ denotes the centre of the groupoid,
$$
Z({}^{\VHS}\!H) = \{z \in  \prod_{x \in X_0} {}^{\VHS}\!H(x,x)\,:\, z_x f=fz_y \, \forall f \in {}^{\VHS}\!H(x,y)\}.
$$
\end{proposition}
\begin{proof}
As in  \cite{Simpson} Lemma 5.1, the map
$$
\Aut(G_{\alpha}) \to \Hom(\pi_fX, G_{\alpha})
$$
is a closed immersion of schemes, so the map
$$
(S/\bG_m)^{\delta} \to \Aut(G_{\alpha})
$$
is analytic, hence continuous. Since  $S/\bG_m$ is isomorphic to $U(1)$, this means that it defines a one-parameter subgroup, so is algebraic.
Therefore the map
$$
S/\bG_m \by  {}^{\VHS}\!H\to{}^{\VHS}\!H 
$$
is algebraic, as $ {}^{\VHS}\!H=\Lim G_{\alpha}$.

Since $\varpi_f(X)^{\red}$  is equivalent to a disjoint union of reductive proalgebraic groups, $G_{\alpha}$ is equivalent to a disjoint union of reductive algebraic groups.
This implies that the connected component $\Aut(G_{\alpha})^0$ of the identity in $\Aut(G_{\alpha})$ is given by
$$
\Aut(G_{\alpha})^0= \prod_{x \in X_0} G_{\alpha}(x,x)/Z(G_{\alpha}).
$$

Since 
$$
\prod_{x \in X_0}{}^{\VHS}\!H(x,x)/Z({}^{\VHS}\!H)=\Lim \prod_{x \in X_0}G_{\alpha}(x,x)/Z(G_{\alpha}),
$$ 
we have an algebraic  map 
$$
S/\bG_m \to\prod_{x \in X_0}{}^{\VHS}\!H(x,x)/Z({}^{\VHS}\!H),
$$
as required.
\end{proof}

\begin{proposition}\label{vhsequiv}
The following conditions are equivalent:
\begin{enumerate}
\item $V$ is a representation of ${}^{\VHS}\!H$;
\item $V$ is a  representation of $H$ such that $\lambda \dmd V \cong V$ for all $\lambda \in (S/\bG_m)(\R)$;
\item $V$ is a  representation of $H$ such that $\lambda \dmd V \cong V$ for some non-torsion $\lambda \in (S/\bG_m)(\R)$.
\end{enumerate}
\end{proposition}
\begin{proof}
$ $

\begin{enumerate}
\item[1.$\implies$2.] If $V$ is a representation of ${}^{\VHS}\!H$, then it is a representation of $H$, so is a semisimple representation of $\varpi_fX$. By Lemma \ref{vhsalg}, $\lambda \in S/\bG_m$ is an inner automorphism of ${}^{\VHS}\!H$, coming from $g \in \prod_{x \in X_0}{}^{\VHS}\!H(x,x)$, say. Then multiplication by $g$ gives the isomorphism $\lambda \dmd V \cong V$.

\item[2.$\implies$3.] Trivial.

\item[3.$\implies$1.] Let $M$ be the monodromy groupoid of $V$; this is a quotient of  $H$. The isomorphism $\lambda \dmd V \cong V$ gives an element $g \in \Aut(M)$, such that $g$ is the image of $\lambda$  in $\Hom(\pi_fX, M)$, using the standard embedding of  $\Aut(M)$ as a closed subscheme of  $\Hom(\pi_fX, M)$. The same is true of $g^n,\lambda^n$, so  the image of $S/\bG_m$ in $\Hom(\pi_fX, M)$ is just the closure of $\{g^n\}_{n \in \Z}$, which is contained in $\Aut(M)$, as $\Aut(M)$ is closed. For any $\mu \in S/\bG_m$, this gives us an isomorphism $\mu\dmd V \cong V$, as required. 
 \end{enumerate}
\end{proof}

\begin{lemma}
If $\rho: \pi_fX \to H(\R)$ has a real Hodge structure, then  the obstruction to a surjective  map $\alpha: H \to G$ factoring through  ${}^{\VHS}\!H$ lies in $\H^1(X, \rho^*\ad \alpha)$.
\end{lemma}
\begin{proof}
We have a real holomorphic map
$$
S/\bG_m \by \pi_fX \to G,
$$
and $\alpha$ will factor through ${}^{\VHS}\!H$ if and only if the induced map
$$
S/\bG_m \xra{\phi} \Hom(\pi_fX, G)/\Aut(G)
$$ 
is constant. Since $G$ is reductive and $S/\bG_m$ connected, it suffices to replace $\Aut(G)$ by the group of inner automorphisms. On tangent spaces, we then have a map
$$
i\R \xra{D_1\phi } \H^1(X, \rho^*\ad \alpha);
$$ 
let $\varphi \in \H^1(X, \rho^*\ad \alpha)$ be the image of $i$.

If $\phi$ is constant, then $\varphi=0$. Conversely, observe that for $t \in (S/\bG_m)(R)$, $D_t\phi= tD_1\phi t^{-1}$, making use of the discrete action of $(S/\bG_m)(R)$ on $\Hom(\pi_fX, G)$. If $\varphi=0$, this implies that $D_t\phi=0$ for all $t \in (S/\bG_m)(\R)$, so $\phi$ is constant, as required.
\end{proof}

\begin{remark}\label{higgsobs}
For compact K\"ahler manifolds, we will have $\varphi=[i\theta -i\bar{\theta}]$, for $\theta$ the Higgs form.
\end{remark}

\begin{definition}
Define ${}^{\VHS}\!G(X)^{\rho, \mal}$ to be the relative Malcev completion of $G(X) \to {}^{\VHS}\!H $. If $\rho: \pi_f \to \varpi_f(X)^{\red}(\R)$ is the canonical map, we write  ${}^{\VHS}\!G(X)^{\alg}$ for  ${}^{\VHS}\!G(X)^{\rho, \mal}$. 
\end{definition}

\begin{lemma}
Given $\lambda \in S(\R)$ such that $\lambda$ is not torsion in $S/\bG_m$, lift the image of $\lambda$ in $\Out(G(X)^{\rho, \mal})$ to some $\alpha \in \Aut(\m \rtimes H)$, where $(\m \rtimes H)$ is some minimal model for $G(X)^{\rho, \mal}$. Then 
$$
{}^{\VHS}\!G(X)^{\rho, \mal}={}^{\alpha}\!(\m \rtimes H). 
$$
\end{lemma}
\begin{proof}
Combine Proposition \ref{vhsequiv} with \cite{weight1}
Corollary \ref{weight1-algboth} and Lemma \ref{weight1-algchar}.
\end{proof}

\subsection{Classical homotopy groups}

\begin{proposition}\label{hodgeclassicalpi}
If a topological space $X$ has a real Hodge structure and satisfies the conditions of Theorem \ref{classicalpimal} for  $R=(\pi_fX)^{\red}$ (or any quotient to which the $U(1)$-action descends), then the homotopy groups $\pi_n(X,x)$ for $n\ge 2$ carry natural real Hodge structures in the sense of \cite{Hodge2}, i.e. algebraic $\Cx^*$-actions on  the real vector spaces
$$
\pi_n(X,x)\ten_{\Z}\R.
$$
These are unique up to conjugation by $\Ru(\varpi_1(X,x))$ (in other words, by the unipotent radical of the Zariski closure of $\pi_1(X,x) \to \GL(\pi_n(X,x)\ten_{\Z}\R)$).
\end{proposition}
\begin{proof}
The real Hodge structure on $X$ defines a discrete action of $\Cx^*$ on $\varpi_n(X,x) \cong \pi_n(X,x)\ten_{\Z}\R$, for which the action of $\R^*$ is algebraic, and the map
$$
 \Cx^* \by \pi_n(X,x) \to \pi_n(X,x)\ten_{\Z}\R
$$
is real holomorphic. Since $\pi_n(X,x)\ten_{\Z}\R$ is finite-dimensional, the argument of  Proposition \ref{vhsalg} adapts to show that the action of $U(1) < \Cx^*$ must be algebraic. Therefore the action of $\Cx^*= (U(1)\by \R^*)/(-1,-1)$ is algebraic.
\end{proof}

\section{Compact K\"ahler manifolds} \label{cK}

\subsection{Formality}\label{formality}

As in 
 \cite{htpy} Theorem \ref{htpy-kformal}, we know that 
the  real pro-algebraic homotopy type of a compact K\"ahler manifold is formal. We recall the proof here.

Given a semisimple complex local system $\ww$ on a compact K\"ahler manifold $X$, there exists a harmonic metric $K$ on $\ww$ (\cite{Simpson} Theorem 1). In fact, this metric is  unique up to a scalar on complex irreducible local systems. In \cite{Simpson} \S 1, the connection $D$ and metric $K$ are used to define connections $D'_K, D''_K$ on $\sA^0_X(\ww)$. Uniqueness implies that the restrictions of $D'_K, D''_K$ to any irreducible subsystem  are independent of the choice of harmonic metric on $\vv$, so $D'_K, D''_K$ are independent of the choice of harmonic metric on $\ww$, and we will simply denote them by $D',D''$.

Now define an operator  $D^c:=iD' -iD''$ on the bundle $\sA^0_X(\ww)$. This is a real operator; in other words, if $\vv$ is a real semisimple local system, then  $D^c:\vv\ten \Cx \to \sA^1_X(\vv\ten \Cx) $ descends to an $\R$-linear map $D^c:\vv \to \sA^1_X(\vv)$. 

Now, the principle of two types (\cite{Simpson} Lemmas 2.1 and 2.2) gives   quasi-isomorphisms
$$
A^{\bt}(X,\vv) \leftarrow (\z_{D^c}(A^*(X,\vv)),d) \to (\H_{D^c}^*(X,\vv),0)
$$
for any semisimple local system $\vv$. As in Proposition \ref{formalpin}, this gives a natural weight decomposition on $G(X)^{\alg}$,  defined by setting $\H^i(X,\vv)$ to be of weight $i$.

\subsection{$S$-action}\label{saction}

In this section, we will show how to construct the map
$$
S(\R) \to \Out(G(X)^{\alg}),
$$
for $X$ a compact K\"ahler manifold.

To every  complex semisimple local system $\vv$, there  corresponds a Higgs bundle $(\sE,\theta)$ with a harmonic metric $K$ (\cite{Simpson} Theorem 1). Here, $\sE$ is a locally free $\sO_X$-module on $X$, and $\theta: \sE \to \sE\ten_{\sO_X}\Omega_X^1$ is a linear map such that $\theta^2=0$. The conditions for a metric to be harmonic are given in \cite{Simpson} \S 1 or \cite{mochi} \S 2.5.

We will proceed by studying the harmonic forms of \cite{Simpson} \S 1. On the space of forms
$$
A^{pq}(X,\sE):=\Gamma(X, \sA^{pq}_{\Cx}\ten_{\O_X}\sE),
$$
we  have an operator $\bar{\pd}$ of type $(0,1)$ coming from the holomorphic structure of $\sE$. The Higgs form $\theta$ then defines a linear operator of type $(1,0)$ on $A^*(X,\sE)$. 

If we set $D'' =\bar{\pd} + \theta$,
the Laplacian is defined by 
$$
\Delta_{\theta}=D''(D'')^* + (D'')^*D'',
$$ 
where $(D'')^*$ is the formal adjoint of $D''$ with respect to the metric $K$.
The cohomology groups $H^i(X,\vv)$ are isomorphic to the space  
$$
\cH^i(\sE,\theta):=\ker \Delta_{\theta} < A^i(X,\sE)
$$ 
of harmonic forms. 

\begin{definition}
Define an action of the circle group  $U(1)$  on the space of semisimple complex local systems by defining $t\vv$ to be the local system corresponding to the Higgs bundle $(\sE,t\theta)$, where $(\sE,\theta)$ is the Higgs bundle corresponding to $\vv$.
\end{definition}

Since $\Delta_{\theta}$ is not of type $(0,0)$, the space $\cH^i(\sE,\theta)$ will not be preserved by the $S$-action. However, we have the following lemma:

\begin{lemma}\label{dmddelta}
Given $\lambda \in \Cx^*$, and a semisimple complex local system $\vv$ on $X$, 
$$
\lambda\dmd\cH^i(\vv)=\cH^i(\frac{\lambda}{\bar{\lambda}}\vv).
$$
\end{lemma}
\begin{proof}
First take $t \in U(1)$; 
we can then   decompose $\Delta_{\theta}=\Delta_{\theta}^{0,0} +\Delta_{\theta}^{1,-1}+\Delta_{\theta}^{-1,1}$ by type, so that
\begin{eqnarray*}
\Delta_{\theta}^{0,0} &=&\bar{\pd}\bar{\pd}^* +\bar{\pd}^*\bar{\pd} + \theta\theta^* +\theta^*\theta\\
\Delta_{\theta}^{1,-1} &=& \theta\bar{\pd}^* +\bar{\pd}^*\theta\\
\Delta_{\theta}^{-1,1}&=& \bar{\pd}\theta^*+\theta^*\bar{\pd}.
\end{eqnarray*}

As remarked in \cite{Simpson} Lemma 4.4, for $t \in U(1)$, $ (\sE,t\theta)$ is another Higgs bundle for which $K$ is a harmonic metric (note that for general  $t \in \Cx^*$, this is not the case, although there will exist another harmonic metric). We then have
\begin{eqnarray*}
\Delta_{t\theta}^{0,0} &=&\Delta_{\theta}^{0,0}\\
\Delta_{t\theta}^{1,-1} &=& t\Delta_{\theta}^{1,-1}\\
\Delta_{t\theta}^{-1,1}&=& t^{-1}\Delta_{\theta}^{-1,1}
\end{eqnarray*}

However, if $\lambda \in \Cx^*$, then
$$
\lambda \dmd \Delta_{\theta}(\lambda^{-1}\dmd v)= \Delta_{\theta}^{0,0} + t\Delta_{\theta}^{1,-1}+ t^{-1}\Delta_{\theta}^{-1,1},
$$
for $t=\frac{\lambda}{\bar{\lambda}}$.

Therefore 
$$
\lambda\dmd\cH^i(\sE,\theta)=\cH^i(\sE, \frac{\lambda}{\bar{\lambda}}\theta),
$$
as required.
\end{proof}

\begin{proposition}\label{dmdcoho}
There is a  family of isomorphisms
$$
\H^i(X,\vv)\xra{\lambda\dmd} \H^i(X, \frac{\lambda}{\bar{\lambda}}(\vv)),
$$
 pure of weight $i$, parametrised by $\lambda \in S(\R) \cong\Cx^*$. This respects cup-products, in the sense that
$$
(\lambda\dmd a)\cup (\lambda \dmd b) = \lambda\dmd (a\cup b) \in \H^*(X, \frac{\lambda}{\bar{\lambda}}(\vv\ten \ww)).
$$ 
\end{proposition}
\begin{proof}
This is defined by using the isomorphism between harmonic forms and cohomology.
 Observe that if $\lambda$ is real, then this action is just multiplication by $\lambda^i$, so is  pure.   

For the last part, given $a \in \cH^i(X, \vv),\, b \in \cH^j(X, \ww)$, let $c \in \cH^{i+j}(X, \vv \ten \ww)$ be uniquely determined by the property that $[c]=[a \wedge b] \in \H^{i+j}(X, \vv \ten \ww)$. Then the principle of two types (\cite{Simpson} Lemma 2.1) implies that $a\wedge b -c \in \im D'_{\theta}D''_{\theta}$, where $D'_{\theta}=\pd +\bar{\theta}$. Since 
$$
\lambda \dmd (D'_{\theta}D''_{\theta}x) = \lambda \bar{\lambda}D'_{\frac{\lambda}{\bar{\lambda}}\theta}D''_{\frac{\lambda}{\bar{\lambda}}\theta} (\lambda \dmd x) ,
$$ 
we have $[\lambda \dmd D'_{\theta}D''_{\theta}x]=0$, so 
$$
\lambda \dmd ([a] \cup [b])= [\lambda\dmd c] =[\lambda \dmd (a \wedge b)]= [\lambda \dmd a] \cup [\lambda \dmd b].
$$  
\end{proof}

By \cite{Simpson} Lemma 2.11, the complex conjugate bundle to $(\sE,\theta)$ is isomorphic to $(\sE^{\vee}, \theta^t)$, with the isomorphism given by the metric, which means that 
$$
\overline{(t\vv)} \cong t\overline{\vv},
$$
for all $t \in U(1)$.

\begin{proposition}
There is a commutative diagram
$$
\begin{CD}
\H^i(X,\vv)@>{\lambda\dmd}>> \H^i(X, \frac{\lambda}{\bar{\lambda}}(\vv))\\
@VVV @VVV\\
\H^i(X,\overline{\vv}) @>{\lambda\dmd}>> \H^i(X,\frac{\lambda}{\bar{\lambda}}(\overline{\vv})),
\end{CD}
$$
where the vertical arrows combine complex conjugation with the isomorphisms $\overline{(t\vv)} \xra{C_t} t(\overline{\vv})$. In particular, the action of $S(\R)$ defines a real Hodge structure.
\end{proposition}
\begin{proof}
Note that
$$
\overline{\lambda'\bar{\lambda}''v}=\lambda'\bar{\lambda}''\bar{v},
$$
for all forms.
\end{proof}

\begin{definition}\label{sactiondef}
We now define the $S$-action on $\pi_f(X)^{\red}$ by Tannakian duality, setting $\lambda\dmd \vv:=\frac{\lambda}{\bar{\lambda}}\vv$, for a real semisimple $\pi_f(X)^{\red}$-representation $\vv$, noting that this is well-defined by the above discussion.

 This action gives us a canonical isomorphism $\lambda\dmd O(\pi_f(X)^{\red}) \cong O(\pi_f(X)^{\red})$, allowing us to define an $S$ action on $\H^*(X, \bO(\pi_f(X)^{\red}))$, as the composition
$$
\H^*(X, \bO(\pi_f(X)^{\red})) \xra{\lambda \dmd} \H^*(X, \lambda \dmd \bO(\pi_f(X)^{\red})) \cong \H^*(X, \bO(\pi_f(X)^{\red})).
$$
\end{definition}

We have therefore defined an action  $S(\R)\to \Out(G(X)^{\alg})$, by universality. We already know (from \cite{htpy} \S \ref{htpy-hodge}) that $\bG_m < S$ acts algebraically. It only remains to show holomorphy. 

\begin{remark}
In \cite{KTP}, a discrete $\Cx^*$-action was defined on the complex schematic homotopy type $G(X)^{\alg}\ten_{\R}\Cx$ of $X$, which leads to the question of whether that action can be compared  to the $S(\R)\cong \Cx^*$-action given here on the real schematic homotopy type. On cohomology groups, the answer is that, for $\lambda \in U(1) < \Cx^*$, the action defined in \cite{KTP} corresponds to $\lambda''$ in our notation. 

Just by considering $\varpi_f(X)^{\red}$, we can see that rest of the action is incomparable with ours in general. However, on the quotient ${}^{\VHS}\!G(X)^{\alg}\ten_{\R}\Cx$ of $G(X)^{\alg}\ten_{\R}\Cx$, our $S$-action extends (by algebraicity) to an action of $S(\Cx) = \Cx^*\by \Cx^*$. For $\lambda \in \Cx^*$, the action defined in \cite{KTP} there corresponds on cohomology groups to $(1,\lambda) \in \Cx^*\by \Cx^*$. 
\end{remark}

\subsection{Holomorphy of the $S$-action}

We now construct the whole Hodge structure, in the form of an element of $O\Hom(X, \bar{W} G(X)^{\alg})(S_{\hol}) $. Initially, this entails defining an element of $\Hom(X, B\varpi_f^{\red}(X))((S/\bG_m)_{\hol})$, for which we adapt the    approach of \cite{Sim2} \S 7.

\subsubsection{Action on the reductive fundamental groupoid}
Let $R:=\varpi_f(X)^{\red}$, and let  $\bB$  be the universal  $R$-torsor   on $X$, coming from the canonical element of $\Hom(\pi_fX, R(\R))$. By Tannakian duality, the adjoint bundle $\ad \bB$ is given by
$$
\ad \bB = \{ \alpha \in \prod_{\vv} \End(\vv)\,:\, \alpha_{\vv\ten \ww}=\alpha_{\vv}\ten \id \oplus \id \ten \alpha_{\ww}, \, \alpha_{\vv \oplus \ww} =\alpha_{\vv} \oplus \alpha_{\ww}, \, f\alpha_{\vv} =\alpha_{\ww}f\},
$$
where $\vv$ ranges over all semisimple local systems on $X$, and $f$ over all elements of $\Hom(\vv,\ww)$. The Higgs forms thus combine to give
$$
\theta+\bar{\theta} \in \ad \bB\ten \sA^1_X,
$$
and we also have a connection
$$
\pd +\bar{\pd}: \bB(\sA^0_X) \to \ad \bB\ten \sA^1_X.
$$

Let
\begin{eqnarray*}
&\sA^0_X\sO_{S/\bG_m}\ten\Cx:= \\ &\ker(\bar{\pd}_{(S/\bG_m)(\Cx)}: \sA_{X\by (S/\bG_m)(\Cx)}^0\ten \Cx \to \sA_{X\by (S/\bG_m)(\Cx)}^0\ten_{\sA^0_{(S/\bG_m)(\Cx)}}\sA^1_{(S/\bG_m)(\Cx)}\ten \Cx),\\
&\sA^n_X\sO_{S/\bG_m}\ten\Cx:= (\sA_X^0\sO_{S/\bG_m}\ten\Cx)\ten_{\sA^0_X}\sA^n_X,
\end{eqnarray*}
the sheaf of those smooth functions on $X\by (S/\bG_m)(\Cx)$ which are holomorphic along $S/\bG_m$. Write $t \in O(S/\bG_m)^{\hol}\ten \Cx$ for the canonical isomorphism $t: (S/\bG_m)(\Cx) \to \Cx^*$, given by $(a,b)\mapsto \frac{a+ib}{a-ib}$ in the notation of \S \ref{RHS}.

We now define a connection 
$$
D_h:\bB(\sA^0_X\sO_{S/\bG_m}\ten\Cx) \to \ad \bB \ten  \sA^1_X\sO_{S/\bG_m}\ten\Cx\\
$$
on the $R(\sA^0_X\sO_{S/\bG_m}\ten\Cx)$-torsor $\bB(\sA^0_X\sO_{S/\bG_m}\ten\Cx)$
by 
$$
D_{h}=\pd +\bar{\pd} +t\theta +t^{-1}\bar{\theta}.
$$ 
$\ker D_{h}$ is then an $R(\sO_{(S/\bG_m)}^{\hol}\ten \Cx)$-torsor on $X \by {(S/\bG_m)(\Cx)}$, with the isomorphisms $\bB_x \cong R(x,x)$ for each $x \in X$ giving isomorphisms $(\ker D_{h})_x \cong R(\sO_{S/\bG_m}^{\hol}\ten \Cx)(x,x)$ on ${(S/\bG_m)(\Cx)}$. 

If $\sigma$ denotes the complex conjugation map on $(S/\bG_m)(\Cx)$, corresponding to $t\mapsto \bar{t}^{-1}$, then there is a canonical isomorphism $\overline{\ker D_{h}} \cong \sigma_*\ker D_{h}$. Using this isomorphism, let
$$
\bB_{h}:= \{ b \in \pi_{X*}\ker D_{h} \,:\, b =\bar{b}\}.
$$ 
This is then an $R((S/\bG_m)^{\hol})$-torsor on $X$, giving 
$$
\rho_{h} \in \Hom(\pi_fX, R(O(S/\bG_m)^{\hol})).
$$

\begin{remark}
Observe that the $\Cx^*$-action of \cite{Simpson} does not  extend this to a $\Gal(\Cx/\R)$-equivariant  ${S/\bG_m}(\Cx)\cong  \Cx^*$-action  on $\varpi_f^{\Cx,\red} X$, since, for $\lambda \in \Cx^*$, $\bar{\lambda(\vv)}=\lambda\bar{\vv}$, whereas if the action were a ${S/\bG_m}(\Cx)$-action, it would satisfy  $\overline{(\lambda\vv)}=\bar{\lambda}^{-1}\bar{\vv}$. However, there is a holomorphic map $S(\Cx) \to \End(\varpi_f^{\Cx} X) $, but this is generally not multiplicative.
\end{remark}

\subsubsection{The full action}\label{fullaction}

Writing $G(X)^{\alg} =R\ltimes U $, it only remains to define an element of 
$$
O\Hom(X, \bar{W}G(X)^{\alg})(S_{\hol}) = \mc(X, (R\ltimes U)(O(S)^{\hol}) )/\Gg(X, U(O(S)^{\hol}).
$$
We already have a canonical element $\omega$ of $\mc(X, R\ltimes U)$, and the argument above defines an element $\omega_h^{\red}$ of 
$$
\mc(X, R(O(S/\bG_m)^{\hol})).
$$ 

 The fibre of 
$$
O\Hom(X, \bar{W}G(X)^{\alg})(S_{\hol}) \to \Hom(X, \bar{W}R)(S_{\hol})
$$
over $\omega_h^{\red}$ is  
$$
\pi( \CC^{\bt}(X, \bO(\bB_h))\ten_{O(S/\bG_m)^{\hol}}O(S)^{\hol} , U),
$$
in the notation of Definition \ref{dgdefgauge}, where we write
$$
\bO(\bB_h):= (O(R)\ten O(S/\bG_m)^{\hol}) \by^{R((S/\bG_m)_{\hol})}\bB_h.
$$

\begin{definition}
Define the complex $\sA^{\bt}_X(\bO(\bB)\ten O(S/\bG_m)^{\hol})$ of sheaves on $X$ by 
$$
 \sA^n_X(\bO(\bB)\ten O(S/\bG_m)^{\hol}):= \{ v \in \pi_{X*}( \sA^n_X\sO_{S/\bG_m}\ten\Cx \ten \bO(\bB))\,:\, \bar{v}=v\},
$$
with differential given by $D=\pd+\bar{\pd}+\theta +\bar{\theta}$. Define the complex $\sA^{\bt}_X(\bO(\bB_h))$ by
$$
 \sA^n_X(\bO(\bB_h)):=\sA^n_X(\bO(\bB)\ten O(S/\bG_m)^{\hol}),
$$
with differential given by $D_h:=\pd +\bar{\pd} +t\theta +t^{-1}\bar{\theta} $. Note that $\sA^{\bt}_X(\bO(\bB)\ten O(S/\bG_m)^{\hol})$ is a flabby resolution of $\bO(\bB)\ten O(S/\bG_m)^{\hol}$, while $\sA^{\bt}_X(\bO(\bB_h))$ is a flabby resolution of $\bO(\bB_h)$.

Let
\begin{eqnarray*}
A^{\bt}_X(\bO(\bB)\ten O(S/\bG_m)^{\hol})&:=& \Gamma(X,\sA^{\bt}_X(\bO(\bB)\ten O(S/\bG_m)^{\hol}))\\
 A^{\bt}_X(\bO(\bB_h))&:=& \Gamma(X,\sA^{\bt}_X(\bO(\bB_h))).
\end{eqnarray*}
\end{definition}

\begin{lemma}
$A^0_X(\bO(\bB)\ten O(S/\bG_m)^{\hol})\ten \Cx$ is the space of those functions $f: X \by (S/\bG_m)(\Cx) \to \bB\ten_{\R}\Cx$ over $X$ for which
$$
f(-,s):X \to \bB\ten_{\R}\Cx
$$
is infinitely differentiable for all $s \in (S/\bG_m)(\Cx)$, and for which all the partial derivatives
$$
\frac{\pd^n}{\pd x_{\alpha_1}\pd x_{\alpha_2}\cdots \pd x_{\alpha_n}}f(x,-):(S/\bG_m)(\Cx) \to \bB\ten_{\R}\Cx$$
are complex analytic, for all $x \in X$.
\end{lemma}
\begin{proof}
It is immediate that all the functions in $A^0_X(\bO(\bB)\ten O(S/\bG_m)^{\hol})\ten \Cx$ satisfy these properties. For the converse, it suffices to show that a function satisfying these properties is infinitely differentiable. This follows because complex analyticity implies infinite differentiability, and continuity of partial derivatives implies differentiability. \end{proof}

\begin{lemma}
Given $f \in A^0_X(\bO(\bB)\ten O(S/\bG_m)^{\hol})\ten \Cx$ and $s_0 \in (S/\bG_m)(\Cx)$, there exists $h \in A^0_X(\bO(\bB)\ten O(S/\bG_m)^{\hol})\ten \Cx$ such that
$$
f(x,s) =f(x,s_0) + (t(s) -t(s_0)) \frac{\pd (f\circ t^{-1})}{\pd t(s)}(x,t(s)) + (t(s) -t(s_0))^2h(x,s).
$$ 
\end{lemma}
\begin{proof}
Set 
$$
h(x,s) = \left\{ \begin{matrix}  \frac{f(x,s) -f(x,s_0) - (t(s) -t(s_0)) \frac{\pd (f\circ t^{-1})}{\pd t(s)}(x,s)}{(t(s) -t(s_0))^2 } & s \ne s_0 \\  \half \frac{\pd^2 (f\circ t^{-1})}{\pd t(s)^2}(x,t(s_0)) & s=s_0 \end{matrix} \right. 
$$
It follows from the preceding lemma that this has the required properties.
\end{proof}

\begin{proposition}\label{extendhol}
Given an $\R$-linear operator $F: A^*_X(\bO(\bB)) \to A^*_X(\bO(\bB))$, we may extend $F$ to an $O(S/\bG_m)^{\hol}$-linear operator on $A^*_X(\bO(\bB)\ten O(S/\bG_m)^{\hol})$ by the formula
$$
F(f)(-,s):= F(f(-,s)),
$$
for $s \in (S/\bG_m)(\Cx)$.

Similarly, given an $\R$-bilinear operator 
$
B: A^*_X(\bO(\bB)) \ten A^*_X(\bO(\bB)) \to \R,
$
we may extend $B$ to an  $O(S/\bG_m)^{\hol}$-bilinear operator
$$
A^*_X(\bO(\bB)\ten O(S/\bG_m)^{\hol})\ten A^*_X(\bO(\bB)\ten O(S/\bG_m)^{\hol}) \to O(S)^{\hol},
$$
by the formula
$$
B(f,g)(s):= B( f(-,s), g(-,s)).
$$
\end{proposition}
\begin{proof}
Observe that the previous lemmas imply that $F(f) \in A^*_X(\bO(\bB))$ and $B(f,g) \in O(S)^{\hol}$,  as required.
\end{proof}

\begin{definition}
Let 
$$
\cH^n_X(\bO(\bB)\ten O(S/\bG_m)^{\hol}):= \ker(\Delta_{\theta}: A^n_X(\bO(\bB)\ten O(S/\bG_m)^{\hol})\to A^n_X(\bO(\bB)\ten O(S/\bG_m)^{\hol}),
$$
and let
$$
\cH^n_X(\bO(\bB_h)):= \ker(\Delta_{t\theta}:A^n_X(\bO(\bB_h))\to A^n_X(\bO(\bB_h))).
$$
\end{definition}

\begin{proposition}\label{sumharm}
There is a direct sum decomposition
$$
A^n_X(\bO(\bB)\ten O(S/\bG_m)^{\hol})=\cH^n_X(\bO(\bB)\ten O(S/\bG_m)^{\hol}) \oplus \Delta_{\theta}A^n_X(\bO(\bB)\ten O(S/\bG_m)^{\hol}).
$$
\end{proposition}
\begin{proof}
We extend the standard inner product on $A^*_X(\bO(\bB))$ to $A^n_X(\bO(\bB)\ten O(S/\bG_m)^{\hol})$ by Proposition \ref{extendhol}, noting that $\Delta_{\theta}$ is self-adjoint with respect to this product. Thus, if $(\Delta_{\theta})^2f=0$ then
$$
(\Delta_{\theta}f,\Delta_{\theta}f)=  (v,(\Delta_{\theta})^2f)=0,
$$
so $ \Delta_{\theta}f(-,s)=0$ for all $s$, hence $ \Delta_{\theta}f=0$, which shows that the two summands have zero intersection. The Green's operator $G$ also extends by Proposition \ref{extendhol}, and $f -\Delta_{\theta} Gf$ is harmonic, as required. 
\end{proof}

\begin{definition}
Define $\lambda,\bar{\lambda} \in O(S) \ten \Cx$ by $a+ib, a-ib$ respectively, in the notation of \S \ref{RHS}. Note that $t=\frac{\lambda}{\bar{\lambda}}$.
\end{definition}

\begin{definition}
 Define an $O(S/\bG_m)^{\hol}\ten_{O(S/\bG_m)}O(S)$-linear  automorphism 
$$
\lambda \dmd :A^*_X(\bO(\bB)\ten O(S/\bG_m)^{\hol})\ten_{O(S/\bG_m)}O(S) \to A^*_X(\bO(\bB)\ten O(S/\bG_m)^{\hol})\ten_{O(S/\bG_m)}O(S)
$$
by the formula of Definition \ref{dmd}.
\end{definition}

\begin{corollary}\label{decomp}
There is a direct sum decomposition
$$
A^n_X(\bO(\bB_h))\ten_{O(S/\bG_m)}O(S) =\cH^n_X(\bO(\bB_h))\ten_{O(S/\bG_m)}O(S) \oplus \Delta_{t\theta}A^n_X(\bO(\bB_h))\ten_{O(S/\bG_m)}O(S).
$$
\end{corollary}
\begin{proof}
The key observation is that  $\lambda \dmd \Delta_{\theta}f = \Delta_{t\theta}(\lambda \dmd f)$, with the same calculation as in Lemma \ref{dmddelta}, so conjugating by $\lambda\dmd$ sends the decomposition of Proposition \ref{sumharm} to this decomposition.
\end{proof}

\begin{proposition}\label{formalhol}
The cosimplicial $R$-representation $\CC^{\bt}(X, \bO(\bB_h)\ten_{O(S/\bG_m)}O(S)) $ in algebras is formal.
\end{proposition}
\begin{proof}
As $\sA^{\bt}_X(\bO(\bB_h)\ten_{O(S/\bG_m)}O(S)$ is a flabby resolution of $\bO(\bB_h)\ten_{O(S/\bG_m)}O(S)$, the argument of Proposition \ref{propforms} shows that $\CC^{\bt}(X, \bO(\bB_h)\ten_{O(S/\bG_m)}O(S)) $ is weakly equivalent to $DA^{\bt}_X(\bO(\bB_h))\ten_{O(S/\bG_m)}O(S)$.

Corollary \ref{decomp} then implies that the $DD^c$ lemma holds for this complex, giving quasi-isomorphisms
$$
\xymatrix@R=-2ex{
(\H^*_{D_{h}^c}(A^{\bt}_X(\bO(\bB_h))\ten_{O(S/\bG_m)}O(S)),0)  \\
&\z_{D_{h}^c}(A^{\bt}_X(\bO(\bB_h))\ten_{O(S/\bG_m)}O(S)),D_h),\ar[lu] \ar[ld]\\ A^{\bt}_X(\bO(\bB_h))\ten_{O(S/\bG_m)}O(S) 
}$$ 
where $D_h^c$ is the conjugate of $D_h$ under the action of $i\dmd$, and $\z_{D_{h}^c}:=\ker(D_{h}^c)$.
\end{proof}

\begin{corollary}
The fibre of 
$$
O\Hom(X, \bar{W}G(X)^{\alg})(S) \to \Hom(X, \bar{W}R)(S)
$$
over $\omega_h^{\red}$ is isomorphic to
$$
\pi(\H^*(X, \bO(\bB_h))\ten_{O(S/\bG_m)^{\hol}}O(S)^{\hol}, NU).
$$ 
\end{corollary}

\begin{definition}
Using the canonical isomorphisms 
$$
\H^*(X, \bO(\bB))\cong \cH^*_X(\bO(\bB)), \quad\H^*(X, \bO(\bB)_h)\ten_{O(S/\bG_m)}O(S)\cong \cH^*_X(\bO(\bB_h))\ten_{O(S/\bG_m)}O(S),
$$ 
define
$$
\lambda \dmd: \H^*(X, \bO(\bB)) \to \H^*(X, \bO(\bB)_h)\ten_{O(S/\bG_m)}O(S),
$$
noting that the argument of  Proposition \ref{dmdcoho} shows that this preserves cup products.
\end{definition}

\begin{definition}\label{hodgedefK}
Recall that $\omega \in  \mc(X, R\ltimes U)  $ corresponds to the canonical map $X \to \bar{W}G(X)^{\alg}$. By \S \ref{formality}, this  gives
$$
\omega \in \pi( \H^*(X, \bO(\bB)), NU),
$$
over $\omega^{\red} $.

Now define 
$$
\omega_h := \lambda \dmd \omega
$$
in 
$$
\pi(\H^*(X, \bO(\bB_h))\ten_{O(S/\bG_m)^{\hol}}O(S)^{\hol}, NU) \cong O\Hom(X, \bar{W}G(X)^{\alg})(S_{\hol})_{\omega_h^{\red}}.
$$
\end{definition}

Thus we have proved:
\begin{theorem}
There is a natural real Hodge structure on any compact K\"ahler manifold $X$.
\end{theorem}
\begin{proof}
Definition \ref{hodgedefK} provides the Hodge structure
$$
h \in O\Hom(X, \bar{W}G(X)^{\alg})(S_{\hol}),
$$
such that $\varpi_f(X)^{\red}$ is pure of weight zero. By \S \ref{saction}, this gives a group homomorphism $S(\R) \to \Out(G(X)^{\alg})$, and by \S \ref{formality}, $\bG_m < S$ gives a weight decomposition.  
\end{proof}

\begin{corollary}
If a compact K\"ahler manifold $X$  satisfies the conditions of Theorem \ref{classicalpimal} for $R=(\pi_fX)^{\red}$ (or any quotient to which the $U(1)$-action descends), then the homotopy groups $\pi_n(X,x)$ for $n\ge 2$ carry natural real Hodge structures in the sense of \cite{Hodge2}, i.e. algebraic $\Cx^*$-actions on  the real vector spaces
$$
\pi_n(X,x)\ten_{\Z}\R.
$$
\end{corollary}
\begin{proof}
Apply Proposition \ref{hodgeclassicalpi}.
\end{proof}

\section{Proper complex varieties}\label{propcx}

In this section, we will show how the techniques of cohomological descent allow us to extend real Hodge structures to all proper complex varieties. By \cite{SD} Remark 4.1.10, the method of \cite{effective} \S 9 shows that a surjective proper morphism of topological spaces is universally of effective cohomological descent. 

\begin{lemma}
If $f:X \to Y$ is a map of compactly generated Hausdorff topological spaces inducing an equivalence on fundamental groupoids, such that $\RR^if_*\vv=0$ for all local systems $\vv$ on $X$ and all $i>0$, then $f$ is a weak equivalence. 
\end{lemma}
\begin{proof}
Without loss of generality, we may assume that $X$ and $Y$ are path-connected. If $\tilde{X}\xra{\pi} X,\tilde{Y}\xra{\pi'} Y$ are the universal covering spaces of $X,Y$, then it will suffice to show that $\tilde{f}:\tilde{X}\to\tilde{Y}$ is a weak equivalence, since the fundamental groups are isomorphic.

As $\tilde{X},\tilde{Y}$ are simply connected, it suffices to show that $\RR^i\tilde{f}_*\Z=0$ for all $i>0$. By the Leray-Serre spectral sequence, $\RR^i\pi_*\Z=0$ for all $i>0$, and similarly for $Y$. The result now follows from the observation that $\pi_*\Z$ is a local system on $X$.
\end{proof}

\begin{proposition}\label{effectivewks}
If $a:X_{\bt} \to S$ is a morphism (of simplicial topological spaces) of effective cohomological descent, then $|a|:|X_{\bt}| \to S$ is a weak equivalence, where $|X_{\bt}|$ is the geometric realisation of $X_{\bt}$.
\end{proposition}
\begin{proof}
We begin by showing that the fundamental groupoids are equivalent. Since $\H^0(|X_{\bt}|,\Z) \cong \H^0(S,\Z)$, we know that $\pi_0|X_{\bt}|\cong \pi_0S$, so we may assume that $|X_{\bt}|$ and $S$ are both connected. 

Now the fundamental group of $|X_{\bt}|$ is isomorphic to the fundamental group of the simplicial set $d\Sing(X_{\bt})$ (the diagonal of the bisimplicial complex given by the singular sets of the  $X_n$). For any group $G$,  the groupoid of $G$-torsors on $|X_{\bt}|$ is thus equivalent to the groupoid of pairs $(T,\omega)$, where $T$ is a $G$-torsor on $X_0$, and the descent datum $\omega:\pd_0^{-1}T \to \pd_1^{-1}T$ is a morphism of $G$-torsors satisfying 
$$
\pd_2^{-1}\omega \circ \pd_0^{-1}\omega=\pd_1^{-1}\omega, \quad \sigma_0^{-1}\omega = 1.
$$
Since $a$ is effective, this groupoid is equivalent to the groupoid of $G$-torsors on $S$, so the fundamental groups are isomorphic.

Given a local system $\vv$ on $|X_{\bt}|$, there is a corresponding $\GL(V)$-torsor $T$, which therefore descends to $S$. Since $\vv=T\by^{\GL(V)}V$ and $T=a^{-1}a_*T$, we can deduce that $\vv=a^{-1}a_*\vv$, so $\RR^ia_*\vv=0$ for all $i>0$, as required.
\end{proof}

\begin{corollary}\label{smoothresn}
Given a proper complex variety $X$, there exists a smooth proper simplicial variety $X_{\bt}$, unique up to homotopy, and a map $a:X_{\bt} \to X$, such that $|X_{\bt}| \to X$ is a weak equivalence. 

If $Y$ is any complex variety, there exists  a smooth proper simplicial variety $X_{\bt}$, a simplicial divisor $D_{\bt} \subset X_{\bt}$ with normal crossings, and a map $(X_{\bt} - D_{\bt}) \to Y$ such that $|X_{\bt} - D_{\bt}| \to Y$ is a weak equivalence.

In fact, we may take each $X_n$ to be projective, and these resolutions are unique up to homotopy.
\end{corollary}
\begin{proof}
Apply \cite{Hodge3} 6.2.8, 6.4.4 and \S 8.2.
\end{proof}

\subsection{Semisimple local systems}

In this section, we will define the real holomorphic $U(1)$-action on a suitable quotient of the real reductive pro-algebraic fundamental groupoid $\varpi_f(X)^{\red}$ of any proper complex variety (or, indeed, of any simplicial proper complex variety). 

Recall that a local system on a simplicial complex $X_{\bt}$ of topological spaces  is  equivalent to the category of pairs $(\vv, \alpha)$, where $\vv$ is a local system on $X_0$, and $\alpha:\pd_0^{-1}\vv \to \pd_1^{-1}\vv$ is an isomorphism of local systems  satisfying 
$$
\pd_2^{-1}\alpha \circ \pd_0^{-1}\alpha=\pd_1^{-1}\alpha, \quad \sigma_0^{-1}\alpha = 1.
$$

\begin{definition}
Given a  simplicial complex $X_{\bt}$ of smooth proper varieties, define $\varpi_f(X_{\bt})^{\red,\norm}$ to be the quotient of $\varpi_f(X_{\bt})^{\red}$ by the image of $\Ru(\pi_f(X_0))$. Its representations consist of normally semisimple local systems on $X_{\bt}$, i.e. semisimple local systems $\ww$ for which $a_0^{-1}\ww$ is also semisimple, for $a_0:X_0 \to X_{\bt}$. 
\end{definition}

\begin{lemma}\label{normwell}
If $f\co X_{\bt} \to Y_{\bt}$ is a homotopy equivalence of  simplicial smooth proper varieties, then   $\varpi_f(|X_{\bt}|)^{\red,\norm}\cong \varpi_f(|Y_{\bt}|)^{\red,\norm}$. 
\end{lemma}
\begin{proof}
Without loss of generality, we may assume that the matching maps  of $f$ are faithfully flat and proper. Topological and algebraic effective descent then imply that $f^{-1}$ induces an equivalence on the categories of local systems, and that $f^*$ induces an equivalence on the categories of quasi-coherent sheaves, and hence on the categories of Higgs bundles. Since   representations of $\varpi_f(|X_{\bt}|)^{\red,\norm}$ correspond to semisimple objects in the category of Higgs bundles on $X_{\bt}$, this completes the proof.
\end{proof}

\begin{definition}
If $X_{\bt} \to X$ is a resolution as in Corollary \ref{smoothresn}, we therefore denote the corresponding reductive algebraic groupoid by $ \varpi_f(X)^{\red,\norm}:=\varpi_f(|X_{\bt}|)^{\red,\norm}$.
\end{definition}

\begin{proposition}\label{norm}
If $X$ is a proper complex variety with a smooth proper resolution $a:X_{\bt} \to X$, then normally semisimple local systems on $X_{\bt}$ correspond to semisimple local systems on $X$ which remain semisimple on pulling back to the normalisation $\pi: X^{\norm}\to X$ of $X$. 
\end{proposition}
\begin{proof}
First observe that  $\varpi_f(|X_{\bt}|)^{\red,\norm}= \varpi_f(X)^{\red}/\langle a_0\Ru(\varpi_f(X_0)) \rangle$. Lemma \ref{normwell} ensures that  $\varpi_f(X_{\bt})^{\red,\norm}$ is independent of the choice of resolution $X_{\bt}$ of $X$, so can be defined as $\varpi_f(X)^{\red}/\langle f\Ru(\varpi_f(Y)) \rangle$ for any smooth projective variety $Y$ and proper faithfully flat $f$.

Now, since $X^{\norm}$ is normal,  we may make use of an observation on \cite{Am} pp.9--10 (due to M. Ramachandran). $X^{\norm}$ has a proper faithfully flat morphism $g$ from a smooth variety $Y$ with connected fibres over $X^{\norm}$, so the map $\pi_fg\co \pi_fY \to \pi_fX^{\norm}$ is full (from the long exact sequence of homotopy). Thus $\pi_fg$ is surjective, so $g(\Ru  (\varpi_f(Y)) = \Ru \varpi_f(X^{\norm})$. 

Taking $f:Y \to X$ to be the composition $Y \xra{g} X^{\norm} \xra{\pi} X$, we see that  $f\Ru(\pi_f(Y))= \pi\Ru(\pi_f(X^{\norm}))$, so
$\varpi_f(X)^{\red,\norm}= \varpi_f(X)^{\red}/\langle \pi\Ru(\pi_f(X^{\norm})) \rangle$, as required.
\end{proof}

\begin{proposition}\label{properred}
If $X$ is a proper complex variety, then there is an action of the circle group $U(1)$ on $\varpi_f(X)^{\red,\norm}$, such that the composition $U(1) \by \pi_fX \to \varpi_f(X)^{\red,\norm}$ is real holomorphic.
\end{proposition}
\begin{proof}
The key observation is that the $U(1)$-action defined in \cite{Simpson} is functorial in $X$, and that semisimplicity is preserved by pullbacks between smooth proper varieties (since Higgs bundles pull back to Higgs bundles), 
so there is a canonical isomorphism $t(\pd_i^{-1}\vv)\cong \pd_i^{-1}(t\vv)$; thus it makes sense for us to define
$$
t(\vv,\alpha):= (t\vv, t(\alpha)).
$$
By Tannakian duality, this defines a $U(1)$-action on $\varpi_f(X)^{\red,\norm}$. 

Since $X_0,X_1$ are smooth and proper, the actions of $S/\bG_m \cong U(1)$ on their reductive pro-algebraic fundamental groupoids are real holomorphic, corresponding to maps 
$$
\pi_f(X_i) \to \varpi_f(X_i)^{\red}(O(S/\bG_m)^{\hol}).
$$
The morphisms $\varpi_f(X_i) \to \varpi_f(X)$ then give us maps
$$
\pi_f(X_i) \to \varpi_f(X)^{\red,\norm}(O(S/\bG_m)^{\hol}),
$$
compatible with $\pi_f(\pd_j), \pi_f(\sigma_j)$. Since 
$$
\xymatrix@1{\pi_f(X_1) \ar@<.5ex>[r]^-{\pd_0} \ar@<-.5ex>[r]_-{\pd_1}& \pi_f(X_0) \to \pi_fX}
$$
is a coequaliser diagram in the category of groupoids, this gives us a map
$$
\pi_f(X) \to \varpi_f(X)^{\red,\norm}(O(S/\bG_m)^{\hol}),
$$
as required.
\end{proof}

\begin{remark}\label{simpprop}
Note that this holds for any space $X$ which can be resolved by a simplicial smooth proper variety, so that we may take $X$ to be the realisation of any simplicial proper variety, as in \cite{Hodge3}.  
\end{remark}

\subsection{The $S$-action}
Fix $X$ as in Remark \ref{simpprop}
We have now defined the $S/\bG_m \cong U(1)$-action on $\varpi_f(X)^{\red,\norm}$, and we must now extend it to the whole of $G(X)^{\norm}:=G(X)^{\rho, \mal}$, for $\rho:\pi_f(X) \to  \varpi_f(X)^{\red,\norm}$.

\begin{definition}
Write $G(X)^{\norm} =R \ltimes U$, with $R$ reductive and $U$ unipotent. The representation 
$$
\pi_fX \to R(O(S/\bG_m)^{\hol})
$$
from Corollary \ref{properred} corresponds to a $R(O(S/\bG_m)^{\hol})$-torsor $\bB_h$ on $X$

Similarly, write $\bB$ for the $R(\R)$-torsor on $X$ corresponding to the canonical map $\pi_f \to \varpi_f(X)^{\red,\norm}(\R)$.
\end{definition}
As in \S \ref{fullaction}, we now need to define an element of 
$$
\pi( \CC^{\bt}(X, \bO(\bB_h))\ten_{O(S/\bG_m)^{\hol}}O(S)^{\hol} , U),
$$

Since $X_{\bt} \xra{a} X$ is a cohomological descent morphism, the map
$$
\CC^{\bt}(X, \bO(\bB_h)) \xra{a^{-1}} \diag\CC^{\bt}(X_{\bt}, a^{-1}\bO(\bB_h))
$$
is a quasi-isomorphism in $c\Alg(R)$. 

Now,  $a_n^{-1}\vv$ is semisimple on $X_n$ for any normally semisimple local system $\vv$ on $X$ (being a pullback of $a_0^{-1}\vv$). This means that the quasi-isomorphisms of Proposition \ref{formalhol} apply to give quasi-isomorphisms
$$
\CC^{\bt}(X_n, a_n^{-1}\bO(\bB_h)) \sim D\H^*(X_n, a_n^{-1}\bO(\bB_h))
$$
for all $n$, functorial in $\pd_i, \sigma_i$. These therefore combine to give a quasi-isomorphism
$$
\diag\CC^{\bt}(X_{\bt}, a^{-1}\bO(\bB_h)) \sim \diag D\H^*(X_{\bt}, a^{-1}\bO(\bB_h)).
$$

Similarly, there are quasi-isomorphisms
$$
\CC^{\bt}(X, \bO(\bB)) \xra{a^{-1}}\diag\CC^{\bt}(X_{\bt}, a^{-1}\bO(\bB)) \sim \diag D\H^*(X_{\bt}, a^{-1}\bO(\bB)).
$$

\begin{theorem}
If $X$ is a proper complex variety, or the realisation of a simplicial proper complex variety, then $X$ admits a canonical real Hodge structure over $\varpi_f(X)^{\red,\norm}$.
\end{theorem}
\begin{proof}
For the weight decomposition, we just use the quasi-isomorphism
$$
\CC^{\bt}(X, \bO(\bB)) \sim \diag D\H^*(X_{\bt}, a^{-1}\bO(\bB)),
$$
defining $\H^n(X_m)$ to be of weight $n$.

The adjunction map $G(X) \to G(X)^{\alg}(\R) $ gives rise to an element 
$$
\omega \in \pi( \CC^{\bt}(X, \bO(\bB)), U).
$$
This is equivalent to giving an element of 
$$
\pi(\diag D\H^*(X_{\bt}, a^{-1}\bO(\bB)), U),
$$
and we use the $\dmd$-action of \S \ref{fullaction} to set
\begin{eqnarray*}
\omega_h := \lambda \dmd \omega &\in& \pi(\diag D\H^*(X_{\bt}, a^{-1}\bO(\bB_h))\ten_{O(S/\bG_m)^{\hol}}O(S)^{\hol}, U) \\
&\cong& \pi( \CC^{\bt}(X, \bO(\bB_h))\ten_{O(S/\bG_m)^{\hol}}O(S)^{\hol}, U).
\end{eqnarray*}
All of the necessary properties now follow from the corresponding results in \S \ref{cK}.
\end{proof}

\begin{corollary}
If a proper complex variety $X$  satisfies the conditions of Theorem \ref{classicalpimal} for  $R=(\pi_fX)^{\red,\norm}$ (or any quotient to which the $U(1)$-action descends), then the homotopy groups $\pi_n(X,x)$ for $n\ge 2$ carry natural real Hodge structures in the sense of \cite{Hodge2}, i.e. algebraic $\Cx^*$-actions on  the real vector spaces
$$
\pi_n(X,x)\ten_{\Z}\R,
$$
unique up to  automorphism by $\Ru\varpi_1(X,x)$.
\end{corollary}
\begin{proof}
Apply Proposition \ref{hodgeclassicalpi}.
\end{proof}

\bibliographystyle{alphanum}
\addcontentsline{toc}{section}{Bibliography}
\bibliography{references.bib}
\end{document}